\documentclass[11pt,a4paper]{article}
\usepackage{amssymb}
\usepackage{amsmath}
\usepackage{latexsym}
\usepackage{graphicx} 
\usepackage{hyperref}
\usepackage{amsthm}
\usepackage{verbatim}
\usepackage{psfrag}
\usepackage{bbm}
\usepackage{dsfont}
\usepackage{mathrsfs}
\usepackage{tikz,pgfplots}
\usepackage{yfonts}
\usepackage{color}
\usepackage{enumerate}
\usepackage{subcaption}
\usepackage{caption}
\usepackage{pgfplots}
\usepackage{subcaption}
\usepackage{multirow}
\usepackage{titling}

\newcommand*{\collauthor}[2]{{#1}$^{#2}$}

\newcommand*{\affiliation}[2]{$\mbox{}^{{#2}}${#1}}
\newcommand*{\colltitle}[1]{\textbf{#1}}

\newtheorem{lemma}{Lemma}[section]
\newtheorem{theorem}[lemma]{Theorem}

\newenvironment{keywords}[1]{\vspace{1cm}\\{\bf \slshape{Keywords}}\quad\slshape{#1}}{}

\newcommand{\tS}{\tilde{S}}
\newcommand{\ts}{\tilde{s}}

\newcommand{\sr}{\rule[-0.45cm]{0pt}{0.9cm}}

\newenvironment{alphafootnotes}
  {\par\edef\savedfootnotenumber{\number\value{footnote}}
   
   \setcounter{footnote}{0}}
  {\par\setcounter{footnote}{\savedfootnotenumber}}

\begin{document}
\begin{center}
\begin{Large}
  \colltitle{A Multi-step Scheme based on Cubic Spline for solving Backward Stochastic Differential Equations}
\end{Large} 
\vspace*{1.5ex}
\begin{alphafootnotes}
\begin{sc}
\begin{large}
\collauthor{Long Teng}{1,}\footnote{Corresponding author (teng@math.uni-wuppertal.de)},
\collauthor{Aleksandr Lapitckii}{},
\collauthor{Michael G\"UNTHER}{1}
\end{large}
\end{sc}
\end{alphafootnotes}
\vspace{1.5ex}

\affiliation{Lehrstuhl f\"ur Angewandte Mathematik und Numerische Analysis,\\
Fakult\"at f\"ur Mathematik und Naturwissenschaften,\\
Bergische Universit\"at Wuppertal, Gau{\ss}str. 20, 
42119 Wuppertal, Germany}{1} 
\end{center}

%%%%%%%%%%%%%%%%%%%%%%%%%%%%%%%%%%%%%%%%%%%%%%%%%%%%%%%%%%%%%%%%%%%%%%%%%%%%%%%%%%%%%%%%%%%%%%%%%%%%%%%%%%%%%%%%%%%%%%%%%%
\section*{Abstract}
In this work we study a multi-step scheme on time-space grids proposed by W. Zhao et al. \cite{Zhao2010} for solving backward stochastic differential equations,
where Lagrange interpolating polynomials are used to approximate the time-integrands with given values of these integrands at chosen multiple time levels. 
For a better stability and the admission of more time levels we investigate the application of spline instead of Lagrange interpolating polynomials to approximate the time-integrands.
The resulting scheme is a semi-discretization in the time direction involving conditional expectations,
which can be numerically solved by using the Gaussian quadrature rules and polynomial interpolations on the spatial grids.
Several numerical examples including applications in finance are presented to demonstrate the high accuracy and stability of our new multi-step scheme.
\begin{keywords}
backward stochastic differential equations, multi-step scheme, cubic splines, time-space grid, Gauss-Hermite quadrature rule
\end{keywords}

\section{Introduction}
Recently, the forward-backward stochastic differential equation (FBSDE) becomes an important tool for formulating many problems in, e.g., mathematical
finance and stochastic control. The BSDE exhibits usually no analytical solution, see e.g., \cite{Karoui1997a}. Their numerical solutions have
thus been extensively studied by many researchers.
The general form of (decoupled) FBSDEs reads
\begin{equation}\label{eq:decoupledbsde}
\left\{
\begin{array}{l}
\,\,\, dX_t = a(t, X_t)\,dt + b(t, X_t)\,dW_t,\quad X_0=x_0,\\
-dY_t = f(t, X_t, Y_t, Z_t)\,dt - Z_t\,dW_t,\\
\quad Y_T=\xi=g(X_T),
 \end{array}\right.
 \end{equation}
where $X_t, a \in \mathbb{R}^n,$ $b$ is a $n\times d$ matrix, $W_t=(W^1_t,\cdots, W^d_t)^T$ is a $d$-dimensional Brownian motion (all Brownian motions are independent with each other),
$f(t, X_t, Y_t, Z_t): [0, T] \times \mathbb{R}^n \times \mathbb{R}^m \times \mathbb{R}^{m\times d} \to \mathbb{R}^m$
is the driver function and $\xi$ is the square-integrable terminal condition.
We see that the terminal condition $Y_T$ depends on the final value of a forward stochastic differential equation (SDE).

For $a=0~\mbox{and}~b=1,$ namely $X_t=W_t,$ one obtains a backward stochastic differential equation (BSDE) of the form
\begin{equation}\label{eq:bsde}
 \left\{
 \begin{array}{l}
-dY_t = f(t, Y_t, Z_t)\,dt - Z_t\,dW_t,\\%\quad t\in [0, T),\\  
 \quad Y_T = \xi = g(W_T),
 \end{array} \right.
\end{equation}
where $Y_t \in \mathbb{R}^m$ and $f:[0,T]\times\mathbb{R}^m\times\mathbb{R}^{m\times d}\to \mathbb{R}^m.$
In the sequel of this paper, we investigate the numerical scheme for solving \eqref{eq:bsde}. Note that the developed schemes can be applied also for solving \eqref{eq:decoupledbsde},
where the general Markovian diffusion $X_t$ can be approximated, e.g., by using the Euler-Scheme.

The existence and uniqueness of solution of \eqref{eq:bsde} assuming the Lipschitz conditions on $f, a(t,X_t), b(t,X_t)~\mbox{and}~g$ are proven
by Pardoux and Peng \cite{Pardoux1990, Pardoux1992}. The uniqueness of solution is extended under more general assumptions for $f$ in \cite{Lepeltier1997},
but only in the one-dimensional case. 

In recent years, many numerical methods have been proposed for the FBSDEs and BSDEs. Peng \cite{Peng1991} obtained a direct relation between FBSDEs and
partial differential equations (PDEs), see also \cite{Karoui1997b}. Based on this relation, several numerical schemes are proposed, e.g., \cite{Douglas1996, Ma1994, Milstein2006}.
As probabilistic methods, (least-squares) Monte-Carlo approaches are investigated in \cite{Bender2012, Bouchard2004, Gobet2005, Lemor2006, Zhao2006}, and tree-based approaches
in \cite{Crisan2010, Teng2018}. For numerical approximation and analysis we refer to \cite{Bally1997, Bender2008, Ma2009, Ma2005, Zhang2004, Zhao2010}. And many others, e.g.,
some numerical methods for BSDEs applying binomial tree are investigated in \cite{Ma2002}. The approach based on the Fourier method for BSDEs is developed in \cite{Ruijter2015}. 

In \cite{Zhao2010}, a multi-step scheme is achieved by using Lagrange interpolating polynomials.
However, the number of multiple time levels is restricted, the stability condition cannot be satisfied for a high number of time steps. This is actually to be expected due to 
Runge's phenomenon. For this reason,
we study in this work a stable multi-step scheme by using
the cubic spline polynomials, for numerically solving BSDEs on the time-space grids.
More precisely, we use the cubic spline polynomials to approximate the integrands, which are conditional mathematical expectations derived from the original BSDEs.
For this, we need to know values of integrands at multiple time levels, which can be numerically
evaluated, e.g., using the Gauss-Hermite quadrature and polynomial interpolations on the spatial grids.
We will study the convergence and the error estimates for the proposed multi-step
scheme.

In the next section, we start with notation and definitions and derive in Section 3
the reference equations for our multi-step scheme for the BSDEs.
In Section 4, we introduce the multi-step scheme for their discretizations.
Section 5 is devoted to error estimates.
In Section 6, several numerical experiments on different types of (F)BSDEs including financial applications are provided to show the high accuracy and stability.
Finally, Section 7 concludes this work.
%Section 4 is devoted to firstly compare different numerical algorithms, and secondly to recognize the effect of imposing stochastic correlation on implied volatility.
%Finally, Section 5 concludes this work.

\section{Preliminaries}
Throughout the paper, we assume that $(\Omega, \mathcal{F}, P; \{\mathcal{F}_t\}_{0\leq t\leq T})$ is a complete, filtered probability space.
In this space, a standard $d$-dimensional Brownian motion $W_t$ with a finite terminal time $T$ is defined, which generates the filtration $\{\mathcal{F}_t\}_{0\leq t\leq T},$ i.e., 
 $\mathcal{F}_t=\sigma\{X_s, 0\leq s \leq t\}$ for FBSDEs or $\mathcal{F}_t=\sigma\{W_s, 0\leq s \leq t\}$ for BSDEs. And the usual hypotheses should be satisfied. We denote the set of all $\mathcal{F}_t$-adapted and square
integrable processes in $\mathbb{R}^d$ with $L^2=L^2(0,T;\mathbb{R}^d).$
A pair of process $(Y_t,Z_t): [0,T]\times\Omega \to \mathbb{R}^m \times \mathbb{R}^{m\times d}$ is the solution of the BSDE \eqref{eq:bsde}
if it is $\mathcal{F}_t$-adapted and square integrable and satisfies \eqref{eq:bsde} as
\begin{equation}\label{eq:generalBSDE_01}
 Y_t = \xi + \int_{t}^T f(s, Y_s, Z_s)\,ds - \int_t^T Z_s\,dW_s,\quad t \in [0, T],
\end{equation}
where $f(t, Y_s, Z_s): [0, T] \times \mathbb{R}^m \times \mathbb{R}^{m\times d} \to \mathbb{R}^m$ is $\mathcal{F}_t$ adapted,
$g:\mathbb{R}^d \to \mathbb{R}^m.$ As mentioned above, these solutions exist uniquely under Lipschitz conditions. 

Suppose that the terminal value $Y_T$ is of the form $g(W^{t,x}_T),$ where $W^{t,x}_T$ denotes the value of $W_T$ starting
from $x$ at time $t.$
Then the solution $(Y^{t,x}_t, Z^{t,x}_t)$ of BSDEs \eqref{eq:bsde} can be represented
\cite{Karoui1997b, Ma2005, Pardoux1992, Peng1991} as
\begin{equation}\label{eq:pderelation}
 Y^{t,x}_t = u(t, x), \quad Z^{t,x}_t=\nabla u(t, x) \quad \forall t \in [0, T),
\end{equation}
which is the solution of the semilinear parabolic PDE of the form
\begin{equation}
\frac{\partial u}{\partial t} + \frac{1}{2}\sum_i^d \partial^2_{i,i} u + f(t, u, \nabla u)=0
\end{equation}
with the terminal condition $u(T,x)=g(x).$ 
In turn, suppose $(Y, Z)$ is the solution of BSDEs, $u(t, x)=Y^{t,x}_t$ is a viscosity solution to the PDE.
\section{Reference equations for the multi-step scheme}
In this section we drive the reference equations for the multi-step scheme by using the cubic spline polynomials.
\subsection{The one-dimensional reference equations}
We start with the one-dimensional processes, namely $m=n=d=1.$
We introduce the uniform time partition for the time interval $[0, T]$
\begin{equation}\label{eq:timepartition}
\Delta_t = \{t_i | t_i \in [0, T], i=0,1,\cdots,N_T, t_i<t_{i+1}, t_0=0, t_{N_T}=T \}.
 \end{equation}
Let $\Delta t :=h= \frac{T}{N_T}$ be the time step, and thus $t_i=t_0 + i h,~\mbox{for}~i=0,1,\cdots,N_T.$
Then one needs to discretize the backward process \eqref{eq:generalBSDE_01}, namely
\begin{equation}\label{eq:generalBSDE_02}
 Y_t = \xi + \int_{t}^T f(s, \mathbb{V}_s)\,ds - \int_t^T Z_s\,dW_s,
\end{equation}
where $\xi=g(W_T), \mathbb{V}_s=(Y_s, Z_s).$ 
Let $(Y_t, Z_t)$ be the adapted solution of \eqref{eq:generalBSDE_02}, we thus have
\begin{equation}\label{eq:generalBSDE_03}
 Y_i = Y_{i+k} + \int_{t_i}^{t_{i+k}} f(s, \mathbb{V}_s)\,ds - \int_{t_i}^{t_{i+k}} Z_s\,dW_s,\quad t\in[0, T),
\end{equation}
where $1\leq k \leq K_y \leq N_T$ with two given positive integers $k~\mbox{and}~K_y.$
To obtain the adaptability of the solution $(Y_t, Z_t),$ we use conditional expectations $E_i[\cdot](=E[\cdot|\mathcal{F}_{t_i}]).$
We start finding the reference equation for $Y.$
We take the conditional expectations $E_i[\cdot]$ on the both sides of \eqref{eq:generalBSDE_03} to obtain
\begin{equation}\label{eq:disBSDEY_01}
 Y_i = E_i[Y_{i+k}] + \int_{t_i}^{t_{i+k}} E_i[f(s, \mathbb{V}_s)]\,ds.
\end{equation}
We see that the integrand on the right-hand side of \eqref{eq:disBSDEY_01}
is deterministic of time $s.$
%with respect to the filtration $\mathcal{F}_{t_i}.$
When the values of $\mathbb{V}_s,$ $(y_t,z_t)$ are available
on the time levels $t_{i+1},t_{i+2},\cdots,t_{i+K_y},$ an approximation of the integrand in \eqref{eq:disBSDEY_01} can be found.
In this work we choose the cubic spline interpolant $\tS_{K_y,t_i}(s)$ based on the support values $(t_{i+j},E_i[f(t_{i+j}, Y_{i+j}, Z_{i+j})]), j=0,\cdots,K_y,$
namely we have
\begin{equation}\label{eq:disBSDEY_02}
 \int_{t_i}^{t_{i+k}} E_i[f(s, \mathbb{V}_s)]\,ds=\int_{t_i}^{t_{i+k}}\tS_{K_y,t_i}(s)\,ds + R_y^i
\end{equation}
with the residual
\begin{equation}\label{eq:residual_y}
 R_y^i=\int_{t_i}^{t_{i+k}}\left(E_i[f(s, \mathbb{V}_s)] - \tS_{K_y,t_i}(s)\right)\,ds.
\end{equation}
Then we can calculate 
\begin{equation}\label{eq:spline_int_01}
 \int_{t_i}^{t_{i+k}} \tS_{K_y,t_i}(s) \,ds=\int_{t_i}^{t_{i+k}}\sum_{j=0}^{K_y-1} \ts^y_{t_i,j}(s)\,ds=\sum_{j=0}^{K_y-1}\int_{t_i}^{t_{i+k}} \ts^y_{t_i,j}(s)\,ds %\\
\end{equation}
with
\begin{equation}
\ts^y_{t_i,j}(s)=a^y_j+b^y_j(s-t_{i+j})+c^y_j(s-t_{i+j})^2+d^y_j(s-t_{i+j})^3,
 \end{equation}
where $s \in [t_{i+j},t_{i+j+1}],\,j=0,\cdots,K_y-1.$ We straightforwardly calculate
\begin{equation}\label{eq:bsde_integration}
\begin{split}
 \int_{t_i}^{t_{i+k}} \ts^y_{t_i,j}\,ds&=\int_{t_{i+j}}^{t_{i+j+1}} \ts^y_{t_i,j}(s)\,ds\\
 &=a^y_j h+ \frac{b^y_j h^2}{2} + \frac{c^y_j h^3}{3} + \frac{d^y_j h^4}{4}. 
\end{split}
\end{equation}
Note that $j$ satisfying $k-1<j\leq K_y-1$ results an integral with zero value when $k<K_y.$ And the coefficients $a^y_j,b^y_j,c^y_j~\mbox{and}~d^y_j$
are obtained with the support points $(t_{i+j},E_i[f(t_{i+j}, Y_{i+j}, Z_{i+j})]), j=0,\cdots,K_y$ as
\begin{equation}\label{eq:spline_sys}
\begin{cases}
\tS_{K_y,t_i}(t_{i+j})=E_i[f(t_{i+j},Y_{{i+j}},Z_{{i+j}})] \quad &j=0,...,K_y \\
\ts^y_{t_i,j}(t_{i+j}) = \ts^y_{t_i,j+1}(t_{i+j}) \quad  &j=0,1,...,K_y-2\\
\ts^{'y}_{t_i,j}(t_{i+j}) = \ts^{'y}_{t_i,j+1}(t_{i+j}) \quad  &j=0,1,...,K_y-2\\
\ts^{''y}_{t_i,j}(t_{i+j}) = \ts^{''y}_{t_i,j+1}(t_{i+j}) \quad  &j=0,1,...,K_y-2.
\end{cases}
\end{equation}
Obviously, we need two boundary conditions to solve the system above. Since the values of derivatives of $E_i[f(t_{i+j}, Y_{i+j}, Z_{i+j})]$
are unknown, we could thus choose e.g., the natural boundary conditions or Not-a-Knot conditions depending on the value of $K_y.$ 
Combining \eqref{eq:disBSDEY_01}, \eqref{eq:disBSDEY_02}, \eqref{eq:spline_int_01} and \eqref{eq:bsde_integration} we obtain the reference equation for $Y_i$ (based on those support points)
as:
\begin{equation}\label{eq:reference_y_01}
  Y_i = E_i[Y_{i+k}] + \sum_{j=0}^{Ky-1} \left[ a^y_j h+ \frac{b^y_j h^2}{2} + \frac{c^y_j h^3}{3} + \frac{d^y_j h^4}{4}\right]+R_y^i,
\end{equation}
where the coefficients $a_j^y, b_j^y, c_j^y~\mbox{and}~d_j^y$ will be obtained by solving \eqref{eq:spline_sys} together with appropriate boundary conditions and depend
on $Y_i.$ Therefore, \eqref{eq:reference_y_01} is an implicit scheme.
%This will be explicitly analyzed in Section 4.

We now start with the reference equation for $Z.$ By multiplying both sides of the equation \eqref{eq:generalBSDE_03} by $\Delta W_{i+1}:=W_{t_{i+1}} - W_{t_{i}}$
and taking the conditional expectations $E_i[\cdot]$ on both sides of the derived equation we obtain
\begin{equation}\label{eq:disBSDEZ_01}
 -E_i[Y_{i+l}\Delta W_{i+l}] = \int_{t_{i}}^{t_{i+l}}E_i[f(s,\mathbb{V}_s)\Delta W_s]\,ds - \int_{t_{i}}^{t_{i+l}}E_i[Z_s]\,ds,
\end{equation}
where the It\^o isometry and Fubini's theorem are used, $\Delta W_s=W_s-W_{t_i}$ and the given integers $l~\mbox{and}~K_z$ satisfy
$1\leq l \leq K_z.$ Similarly, we derive the reference equation of $Z$ also based on the support points $(t_{i+j}, E_i[f(t_{i+j}, y_{i+j}, z_{i+j})\Delta w_{i+j}])$
and $((t_{i+j}, E_i[z_{i+j}])$, $j=0,\cdots,K_z.$
%However, due to $\Delta W_{t_i}=0$ and $E_i[Z_i]=Z_i$ we can firstly rewrite \eqref{eq:disBSDEZ_01} as
%\begin{equation}\label{eq:disBSDEZ_02}
% Z_i=E_i[Y_{i+l}\Delta W_{i+l}] + \int_{t_{i+1}}^{t_{i+l}}E_i[f(s,\mathbb{V}_s)\Delta W_s]\,ds - \int_{t_{i+1}}^{t_{i+l}}E_i[Z_s]\,ds.
%\end{equation}
Then, we again use the cubic spline polynomials to approximate the time deterministic integers and obtain 
\begin{equation}\label{eq:z_1}
\begin{split}
 \int_{t_{i}}^{t_{i+l}}E_i[f(t_{s}, Y_{s}, Z_{s})\Delta w_{s}]\,ds&=\int_{t_i}^{t_{i+l}}\tS_{K_{z_1},t_i}(s)\,ds + R_{z_1}^i\\
 &=\sum_{j=0}^{K_z-1}\int_{t_{i}}^{t_{i+l}} \ts^{z_1}_{t_i,j}(s)\,ds + R_{z_1}^i
\end{split}
\end{equation}
with
\begin{align}
R_{z_1}^i&=\int_{t_{i}}^{t_{i+l}}\left(E_i[f(t_{s}, Y_{s}, Z_{s})\Delta w_{s}] - \tS_{K_{z_1},t_i}(s)\right)\,ds,\\
\ts^{z_1}_{t_i,j}(s)&=a^{z_1}_j+b^{z_1}_j(s-t_{i+j})+c^{z_1}_j(s-t_{i+j})^2+d^{z_1}_j(s-t_{i+j})^3 \label{eq:zs_1}
\end{align}
$~\mbox{for}~s \in [t_{i+j},t_{i+j+1}],\,j=0,\cdots,K_z-1,$ and 
\begin{equation}\label{eq:z_2}
\begin{split}
  \int_{t_{i}}^{t_{i+l}}E_i[Z_{s}]\,ds&=\int_{t_{i}}^{t_{i+l}}\tS_{K_{z_2},t_i}(s)\,ds + R_{z_2}^i\\
   &=\sum_{j=0}^{K_z-1}\int_{t_{i}}^{t_{i+l}} \ts^{z_2}_{t_i,j}(s)\,ds + R_{z_2}^i
\end{split}
\end{equation}
with
\begin{align}
R_{z_2}^i &=\int_{t_{i}}^{t_{i+l}}\left(E_i[Z_{s}] - \tS_{K_{z_2},t_i}(s)\right)\,ds,\\
\ts^{z_2}_{t_i,j}(s)&=a^{z_2}_j+b^{z_2}_j(s-t_{i+j})+c^{z_2}_j(s-t_{i+j})^2+d^{z_2}_j(s-t_{i+j})^3\label{eq:zs_2}
 \end{align}
$~\mbox{for}~s \in [t_{i+j},t_{i+j+1}],\,j=0,\cdots,K_z-1$ and we let
\begin{equation}\label{eq:residual_z}
 R_{z}^i:=R_{z_1}^i+R_{z_2}^i.
\end{equation}
Furthermore,
using the relation \eqref{eq:pderelation} and integration by parts it can be verified that
\begin{equation}\label{eq:yzrelation}
 E_i[Y_{i+l}\Delta W_{i+l}] =  l h E_i[Z_{i+1}]. 
\end{equation}
Integrating \eqref{eq:zs_1}, \eqref{eq:zs_2} and combining \eqref{eq:disBSDEZ_01}, \eqref{eq:z_1}, \eqref{eq:z_2}
and \eqref{eq:yzrelation} we obtain the reference equation for $Z_i$ as:
\begin{equation}\label{eq:disBSDEZ_05}
\begin{split}
  0=lhE_i[Z_{i+l}] &+ \sum_{j=0}^{Kz-1} \left[ a^{z_1}_j h+ \frac{b^{z_1}_j h^2}{2} + \frac{c^{z_1}_j h^3}{3} + \frac{d^{z_1}_j h^4}{4}\right]\\
  &-\sum_{j=0}^{Kz-1} \left[ a^{z_2}_j h+ \frac{b^{z_2}_j h^2}{2} + \frac{c^{z_2}_j h^3}{3} + \frac{d^{z_2}_j h^4}{4}\right]+R_z^i,
 \end{split}
\end{equation}
where the coefficients $a^{z_1}_j, b^{z_1}_j, c^{z_1}_j, d^{z_1}_j$ are solutions of
\begin{equation}\label{eq:spline_sys_z_1}
\begin{cases}
\tS_{K_z,t_i}(t_{i+j})=E_i[f(t_{i+j}, Y_{i+j}, Z_{i+j})\Delta W_{i+j}] \quad &j=0,...,K_z \\
\ts^{{z_1}}_{t_i,j}(t_{i+j}) = \ts^{{z_1}}_{t_i,j+1}(t_{i+j}) \quad  &j=0,...,K_z-2\\
\ts^{'{z_1}}_{t_i,j}(t_{i+j}) = \ts^{'{z_1}}_{t_i,j+1}(t_{i+j}) \quad  &j=0,...,K_z-2\\
\ts^{''{z_1}}_{t_i,j}(t_{i+j}) = \ts^{''{z_1}}_{t_i,j+1}(t_{i+j}) \quad  &j=0,...,K_z-2
\end{cases}
\end{equation}
with the appropriate boundary conditions,
and the coefficients $a^{z_2}_j, b^{z_2}_j, c^{z_2}_j, d^{z_2}_j$ are solutions of
\begin{equation}\label{eq:spline_sys_z_2}
\begin{cases}
\tS_{K_z,t_i}(t_{i+j})=E_{i}[Z_{{i+j}}] \quad &j=0,...,K_z \\
\ts^{{z_2}}_{t_i,j}(t_{i+j}) = \ts^{z_2}_{t_i,j+1}(t_{i+j}) \quad  &j=0,...,K_z-2\\
\ts^{'{z_2}}_{t_i,j}(t_{i+j}) = \ts^{'{z_2}}_{t_i,j+1}(t_{i+j}) \quad  &j=0,...,K_z-2\\
\ts^{''{z_2}}_{t_i,j}(t_{i+j}) = \ts^{''{z_2}}_{t_i,j+1}(t_{i+j}) \quad  &j=0,...,K_z-2
\end{cases}
\end{equation}
with the appropriate boundary conditions, respectively. 

\subsection{The high-dimensional reference equations}
In this section, we give the reference equations for the high-dimensional case. With the aid of \eqref{eq:reference_y_01} we can straightforwardly
write the reference equation for $y_i$ in component-wise as
\begin{equation}\label{eq:high_y}
  Y^{\tilde{m}}_i = E_i[Y^{\tilde{m}}_{i+k}] + \sum_{j=0}^{Ky-1} \left[ a^{y_j, \tilde{m}} h+ 
  \frac{b^{y_j,\tilde{m}} h^2}{2} + \frac{c^{y_j,\tilde{m}} h^3}{3} + \frac{d^{y_j,\tilde{m}} h^4}{4}\right]+R_y^{i,\tilde{m}},
\end{equation}
with
\begin{equation}\label{eq:spline_sys_h}
\begin{cases}
\tS^{\tilde{m}}_{K_y,t_i}(t_{i+j})=\mathbb{E}_{i}[f^{\tilde{m}}(t_{i+j},Y_{{i+j}},Z_{{i+j}})] \quad &j=0,...,K_y \\
\ts^{y,\tilde{m}}_{t_i,j}(t_{i+j}) = \ts^{y,\tilde{m}}_{t_i,j+1}(t_{i+j}) \quad  &j=0,1,...,K_y-2\\
\ts^{'y,\tilde{m}}_{t_i,j}(t_{i+j}) = \ts^{'y,\tilde{m}}_{t_i,j+1}(t_{i+j}) \quad  &j=0,1,...,K_y-2\\
\ts^{''y,\tilde{m}}_{t_i,j}(t_{i+j}) = \ts^{''y,\tilde{m}}_{t_i,j+1}(t_{i+j}) \quad  &j=0,1,...,K_y-2,
\end{cases}
\end{equation}
where $f^{\tilde{m}}$ is the $\tilde{m}$-th component of the vector $f$ for $\tilde{m}=1,2,\cdots,m.$
The coefficients $a_j^{y,\tilde{m}}, b_j^{y,\tilde{m}}, c_j^{y,\tilde{m}}~\mbox{and}~d_j^{y,\tilde{m}}$
will be obtained by solving the $\tilde{m}$-th system \eqref{eq:spline_sys_h} together with appropriate boundary conditions.
The $\tilde{m}$-th component residual reads
\begin{equation}\label{eq:residuals_y}
 R_y^{i,\tilde{m}}=\int_{t_i}^{t_{i+k}}\left(E_i[f^{\tilde{m}}(s, Y_s, Z_s)] - \tS^{\tilde{m}}_{K_y,t_i}(s)\right)\,ds.
\end{equation}
Similarly, the reference equation for $Z_i$ can be formulated as follows:
\begin{equation}\label{eq:high_z}
\begin{split}
  0=lhE_i[Z^{\tilde{m},\tilde{d}}_{i+l}] &+ \sum_{j=0}^{Kz-1} \left[ a^{z_1,\tilde{m},\tilde{d}}_j h+ \frac{b^{z_1,\tilde{m},\tilde{d}}_j h^2}{2}
  + \frac{c^{z_1,\tilde{m},\tilde{d}}_j h^3}{3} + \frac{d^{z_1,\tilde{m},\tilde{d}}_j h^4}{4}\right]\\
  &-\sum_{j=0}^{Kz-1} \left[ a^{z_2,\tilde{m},\tilde{d}}_j h+ \frac{b^{z_2,\tilde{m},\tilde{d}}_j h^2}{2} + \frac{c^{z_2,\tilde{m},\tilde{d}}_j h^3}{3} + \frac{d^{z_2,\tilde{m},\tilde{d}}_j h^4}{4}\right]
  + R_{z}^{i,{\tilde{m},\tilde{d}}},
 \end{split}
\end{equation}
where the coefficients $a^{z_1,\tilde{m},\tilde{d}}_j, b^{z_1,\tilde{m},\tilde{d}}_j, c^{z_1,\tilde{m},\tilde{d}}_j, d^{z_1,\tilde{m},\tilde{d}}_j$ are solutions of
\begin{equation}
\begin{cases}
\tS^{\tilde{m},\tilde{d}}_{K_z,t_i}(t_{i+j})=E_i[f^{\tilde{m}}(t_{i+j}, Y_{i+j}, Z_{i+j})\Delta W^{\tilde{d}}_{i+j}] \quad &j=0,...,K_z \\
\ts^{{z_1},\tilde{m},\tilde{d}}_{t_i,j}(t_{i+j}) = \ts^{{z_1},\tilde{m},\tilde{d}}_{t_i,j+1}(t_{i+j}) \quad  &j=0,...,K_z-2\\
\ts^{'{z_1},\tilde{m},\tilde{d}}_{t_i,j}(t_{i+j}) = \ts^{'{z_1},\tilde{m},\tilde{d}}_{t_i,j+1}(t_{i+j}) \quad  &j=0,...,K_z-2\\
\ts^{''{z_1},\tilde{m},\tilde{d}}_{t_i,j}(t_{i+j}) = \ts^{''{z_1},\tilde{m},\tilde{d}}_{t_i,j+1}(t_{i+j}) \quad  &j=0,...,K_z-2
\end{cases}
\end{equation}
with the appropriate boundary conditions,
and the coefficients $a^{z_2,\tilde{m},\tilde{d}}_j, b^{z_2,\tilde{m},\tilde{d}}_j, c^{z_2,\tilde{m},\tilde{d}}_j, d^{z_2,\tilde{m},\tilde{d}}_j$ are solutions of
\begin{equation}
\begin{cases}
\tS^{\tilde{m},\tilde{d}}_{K_z,t_i}(t_{i+j})=E_{i}[Z^{\tilde{m},\tilde{d}}_{{i+j}}] \quad &j=0,...,K_z \\
\ts^{{z_2},\tilde{m},\tilde{d}}_{t_i,j}(t_{i+j}) = \ts^{z_2,\tilde{m},\tilde{d}}_{t_i,j+1}(t_{i+j}) \quad  &j=0,...,K_z-2\\
\ts^{'{z_2},\tilde{m},\tilde{d}}_{t_i,j}(t_{i+j}) = \ts^{'{z_2},\tilde{m},\tilde{d}}_{t_i,j+1}(t_{i+j}) \quad  &j=0,...,K_z-2\\
\ts^{''{z_2},\tilde{m},\tilde{d}}_{t_i,j}(t_{i+j}) = \ts^{''{z_2},\tilde{m},\tilde{d}}_{t_i,j+1}(t_{i+j}) \quad  &j=0,...,K_z-2.
\end{cases}
\end{equation}
The corresponding residual reads
\begin{equation}\label{eq:residuals_z}
 R_{z}^{i,{\tilde{m},\tilde{d}}}=R_{z_1}^{i,{\tilde{m},\tilde{d}}}+R_{z_2}^{i,{\tilde{m},\tilde{d}}}
\end{equation}
with
\begin{equation}
R_{z_1}^{i,{\tilde{m},\tilde{d}}}=\int_{t_{i}}^{t_{i+l}}\left(E_i[f^{{\tilde{m}}}(t_{s}, Y_{s}, Z_{s})\Delta W^{{\tilde{d}}}_{s}] - \tS^{{\tilde{m},\tilde{d}}}_{K_{z_1},t_i}(s)\right)\,ds,
\end{equation}
\begin{equation}
R_{z_2}^{i,{\tilde{m},\tilde{d}}} =\int_{t_{i}}^{t_{i+l}}\left(E_i[Z^{\tilde{m},\tilde{d}}_{s}] - \tS^{\tilde{m},\tilde{d}}_{K_{z_2},t_i}(s)\right)\,ds,
\end{equation}
where $\tilde{m} = 1,2,\cdots,m$ and $\tilde{d}=1,2,\cdots,d.$ Note that, by removing superscripts $\tilde{m}$ and $\tilde{d},$
we can write \eqref{eq:high_y} and \eqref{eq:high_z} in matrix form.
%The resulting expression looks like the expression in one-dimensional case. Therefore, in the sequel of this paper we will discuss the discretizations generally with respect to the reference equations without superscripts.  
\subsection{The cubic spline coefficients}
%Since we need $Z_i$ to compute $Y_i,$  we thus have to define the number of time steps as $K=\max\left\{K_y, K_z \right\}.$
As mentioned before, due to the lack of derivative values of the integrands, we should choose some cubic spline which does not need those derivative values.
Furthermore, it will be shown in the next section that \eqref{eq:high_y} is stable for any positive $k$ and $K_y,$ we thus fix $k=K_y.$
However, \eqref{eq:high_z} is only stable for any positive $K_z$ and $l=1.$
Therefore, in the sequel of this paper we fix $k=K_y$ and $l=1.$

For the reference equation \eqref{eq:spline_sys}, we calculate cubic spline coefficients for different values of $K_y$ as follows.
%We will explain our choice for different values of $K$ as follows.
For notational simplicity, we let $g_{i+j}=E_i[f(t_{i+j},Y_{i+j},Z_{i+j})]$ for $j=0,\cdots,K_y.$
\begin{itemize}
 \item $K_y=1:$ there are only two points available. One can just construct a straight line and obtain
 $a^y_0=g_i,b^y_0=\frac{g_{i+1}-g_i}{h},c^y_0=0, d^y_0=0.$
 Now, we can rewrite \eqref{eq:reference_y_01} as
 \begin{align}
  Y_i &= E_i[Y_{i+K_y}] + \frac{h}{2} g_i + \frac{h}{2} g_{i+1}+R_y^i\\
     &:=E_i[Y_{i+K_y}] + h K_y\sum_{j=0}^{Ky} \gamma_{K_y,j}^{K_y} E_i[f(t_{i+j},Y_{i+j},Z_{i+j})]+R_y^i,
 \end{align}
where $\gamma_{K_y,0}^{K_y}=\gamma_{K_y,1}^{K_y}=\frac{1}{2}.$
\item $K_y=2:$ we can already construct e.g., a natural cubic spline based on three points. The corresponding
coefficients can be calculated as follows.\\

For $\ts^{y}_{t_i,0}(s), s\in[t_i,t_{i+1}]: $
\begin{align*}
 a_0 &= g_i, b_0 =-(5g_i-6g_{i+1} + g_{i+2})/4h\\
 c_0 &=0, d_0 =(g_i-2g_{i+1} + g_{i+2})/4h^3
\end{align*}

For $\ts^{y}_{t_i,1}(s), s\in[t_{i+1},t_{i+2}]: $
\begin{align*}
 a_1 &= g_{i+1}, b_1 =-(g_i- g_{i+2})/2h\\
 c_1 &= (3g_i-6g_{i+1} + 3g_{i+2})/4h^2, d_1 =-(g_i-2g_{i+1} + g_{i+2})/4h^3
\end{align*}
Thus, \eqref{eq:reference_y_01} can be rewritten as
 \begin{align}
  Y_i &= E_i[Y_{i+K_y}] + \frac{3h}{8} g_i + \frac{10h}{8} g_{i+1}  + \frac{3h}{8} g_{i+2}+R_y^i\\
     &:=E_i[Y_{i+K_y}] + h K_y\sum_{j=0}^{Ky} \gamma_{K_y,j}^{K_y} E_i[f(t_{i+j},Y_{i+j},Z_{i+j})]+R_y^i,
 \end{align}
where $\gamma_{K_y,0}^{K_y}=\gamma_{K_y,2}^{K_y}=\frac{3}{16}, \gamma_{K_y,1}^{K_y}=\frac{5}{8}.$

Moreover, for the cubic spline we set the second derivatives of cubic interpolants at boundaries to be zero.
Instead of this, one can also choose a second order polynomial for the whole interval, namely $(t_i,t_{i+2}).$
In this way we obtain the polynomial $p_i(s)$ as
\begin{equation}
 g_i(s-t_i) -\left(\frac{3}{2}g_i-2g_{i+1}+ \frac{1}{2}g_{i+2}\right) (s-t_i) /h + \left(\frac{1}{2}g_i-g_{i+1}+ \frac{1}{2}g_{i+2}\right) (s-t_i)^2 /h^2
\end{equation}
and its integration as
\begin{equation}
 \int_{t_i}^{t_{i+2}}p_i(s)ds= h\frac{g_i+4g_{i+1}+g_{i+2}}{3}.
\end{equation}
By using the second order polynomial we rewrite \eqref{eq:reference_y_01} as 
\begin{align}
  Y_i &= E_i[Y_{i+K_y}] + \frac{h}{3} g_i + \frac{4h}{3} g_{i+1}  + \frac{h}{3} g_{i+2}+R_y^i \nonumber\\
     &:=E_i[Y_{i+K_y}] + h K_y\sum_{j=0}^{Ky} \gamma_{K_y,j}^{K_y} E_i[f(t_{i+j},Y_{i+j},Z_{i+j})]+R_y^i,
 \end{align}
where $\gamma_{K_y,0}^{K_y}=\gamma_{K_y,2}^{K_y}=\frac{1}{6}, \gamma_{K_y,1}^{K_y}=\frac{2}{3}.$

\item $K_y=3:$ for $K_y \geq 3$ we will use the Not-a-knot cubic spline and calculate the corresponding
coefficients as follows.\\

For $\ts^{y}_{t_i,0}(s), s\in[t_i,t_{i+1}]: $
\begin{align*}
 a_0 &= g_i, b_0 =-(11g_i-18g_{i+1} + 9g_{i+2} - 2g_{i+3})/6h\\
 c_0 &=(2g_i-5g_{i+1} + 4g_{i+2} - g_{i+3})/2h^2, d_0 =-(g_i-3g_{i+1} + 3g_{i+2}-g_{i+3})/6h^3
\end{align*}

For $\ts^{y}_{t_i,1}(s), s\in[t_{i+1},t_{i+2}]: $
\begin{align*}
 a_1 &= g_{i+1}, b_1 =-(2g_i+3g_{i+1} - 6g_{i+2} + g_{i+3})/6h\\
 c_1 &= (g_i-2g_{i+1} + g_{i+2})/2h^2, d_1 =-(g_i-3g_{i+1} + 3g_{i+2}-g_{i+3})/6h^3
\end{align*}

For $\ts^{y}_{t_i,2}(s), s\in[t_{i+2},t_{i+3}]: $
\begin{align*}
 a_2 &= g_{i+2}, b_2 =(g_i -6g_{i+1} +3g_{i+2} + 2g_{i+3})/6h\\
 c_2 &= (g_i-2g_{i+1} + g_{i+3})/2h^2, d_2 =-(g_i-3g_{i+1}+ 3 g_{i+2} - g_{i+3})/6h^3
\end{align*}
Thus, \eqref{eq:reference_y_01} can be rewritten as
 \begin{align}
  Y_i &= E_i[Y_{i+K_y}] + \frac{3h}{8} g_i + \frac{9h}{8} g_{i+1}  + \frac{9h}{8} g_{i+2}+\frac{3h}{8} g_{i+2}+R_y^i\\
     &:=E_i[Y_{i+K_y}] + h K_y\sum_{j=0}^{Ky} \gamma_{K_y,j}^{K_y} E_i[f(t_{i+j},Y_{i+j},Z_{i+j})]+R_y^i,
 \end{align}
where $\gamma_{K_y,0}^{K_y}=\gamma_{K_y,3}^{K_y}=\frac{1}{8}, \gamma_{K_y,1}^{K_y}=\gamma_{K_y,2}^{K_y}=\frac{3}{8}.$
\end{itemize}
In an analogous way we can also find coefficients for $K_y\geq3,$ and report them for $1\leq K_y \leq 6$ in Table \ref{tab: coeffs_y}.
\begin{table}[h]
	\centering
    \begin{tabular}{|c|c|c|c|c|c|c|c|}
    \hline
    $K_y$ & \multicolumn{7}{|c|}{$\gamma_{K_y,j}^{K_y}$}  \\ \hline
          &$j=0$ & $j=1$ &$j=2$&$j=3$ &$j=4$ &$j=5$ &$j=6$\\ \hline
    1     & $\frac{1}{2}$  & $\frac{1}{2}$ & & & & & \\    \hline
    2 (Second Order Polynomial)    & $\frac{1}{6}$  & $\frac{2}{3}$ &$\frac{1}{6}$ & & & &\\    \hline
    2 (Natural Cubic Spline )     & $\frac{3}{16}$  & $\frac{5}{8}$ &$\frac{3}{16}$ & & & &\\    \hline
    3     & $\frac{1}{8}$  & $\frac{3}{8}$ &$\frac{3}{8}$ &$\frac{1}{8}$ & & & \\    \hline
    4     & $\frac{1}{12}$  & $\frac{1}{3}$ &$\frac{1}{6}$ &$\frac{1}{3}$ &$\frac{1}{12}$ & &\\    \hline
    5     & $\frac{41}{600}$  & $\frac{19}{75}$ &$\frac{107}{600}$ &$\frac{107}{600}$ &$\frac{19}{75}$ &$\frac{41}{600}$ &\\    \hline
    6     & $\frac{19}{336}$  & $\frac{3}{14}$ &$\frac{15}{112}$ &$\frac{4}{21}$ &$\frac{15}{112}$ &$\frac{3}{14}$ &$\frac{19}{336}$\\    \hline
    \end{tabular}
    \caption{The coefficients $[\gamma_{K_y,j}^{K_y}]_{j=0}^{K_y}$ for $K_y=1,2,\cdots,6.$}
    \label{tab: coeffs_y}
\end{table}

We substitute $l=1$ into \eqref{eq:disBSDEZ_05} and thus obtain
\begin{equation}\label{eq:disBSDEZ_06}
\begin{split}
  0=hE_i[Z_{i+1}] &+ \sum_{j=0}^{Kz-1} \left[ a^{z_1}_j h+ \frac{b^{z_1}_j h^2}{2} + \frac{c^{z_1}_j h^3}{3} + \frac{d^{z_1}_j h^4}{4}\right]\\
  &-\sum_{j=0}^{Kz-1} \left[ a^{z_2}_j h+ \frac{b^{z_2}_j h^2}{2} + \frac{c^{z_2}_j h^3}{3} + \frac{d^{z_2}_j h^4}{4}\right]+R_z^i.
 \end{split}
\end{equation}
Note that both the sum terms in the latter equation have the same structure, they will have the same coefficients. 
We use $g_{i+j}$ for $E_i[f(t_{i+j},Y_{i+j},Z_{i+j})\Delta W_{i+j}]$ and $\tilde{g}_{i+j}$ for $E_i[Z_{i+j}]$ for $j=0,1,\cdots,K_z.$
Similar to the way of calculating the coefficients for the reference equation of $Y_i,$ in the following we calculate the coefficients for \eqref{eq:disBSDEZ_06}.
\begin{itemize}
\item $K_z=1:$ we construct straight lines
 $a^{z_1}_0=g_i,b^{z_1}_0=\frac{g_{i+1}-g_i}{h},c^{z_1}_0=0, d^{z_1}_0=0$
 and $a^{z_2}_0=\tilde{g}_i,b^{z_2}_0=\frac{\tilde{g}_{i+1}-\tilde{g}_i}{h},c^{z_2}_0=0, d^{z_2}_0=0$
 Now, we can rewrite \eqref{eq:disBSDEZ_06} as
 \begin{align}
  &0 = hE_i[Z_{i+1}] + \frac{h}{2} g_i + \frac{h}{2} g_{i+1}- \frac{h}{2} \tilde{g}_i - \frac{h}{2} \tilde{g}_{i+1}+R_z^i\\
     &:=hE_i[Z_{i+1}] + h \sum_{j=0}^{Kz} \gamma_{K_z,j}^{1} E_i[f(t_{i+j},Y_{i+j},Z_{i+j})\Delta W_{i+j}]-h \sum_{j=0}^{Kz} \gamma_{K_z,j}^{1} E_i[Z_{i+j}]+R_z^i,
 \end{align}
where $\gamma_{K_z,0}^{1}=\gamma_{K_z,1}^{1}=\frac{1}{2}.$
\item $K_z=2:$ due to $l=1$ we only need to consider the interval $[t_i,t_{i+1}].$ \\
Using natural cubic splines: $\ts^{z_1}_{t_i,0}(s), \ts^{z_2}_{t_i,0}(s), s\in[t_i,t_{i+1}]: $
\begin{align*}
 a^{z_1}_0 &= g_i, b^{z_1}_0 =-(5g_i-6g_{i+1} + g_{i+2})/4h, c^{z_1}_0 =0, d^{z_1}_0 =(g_i-2g_{i+1} + g_{i+2})/4h^3\\
 a^{z_2}_0 &= \tilde{g}_i, b^{z_2}_0 =-(5\tilde{g}_i-6\tilde{g}_{i+1} + \tilde{g}_{i+2})/4h, c^{z_2}_0 =0, d^{z_2}_0 =(\tilde{g}_i-2\tilde{g}_{i+1} + \tilde{g}_{i+2})/4h^3\\
 \end{align*}
Thus, \eqref{eq:disBSDEZ_06} can be rewritten as
 \begin{align}
  &0 = hE_i[Z_{i+1}] + \frac{7h}{16} g_i + \frac{10h}{16} g_{i+1}  - \frac{h}{16} g_{i+2} - (\frac{7h}{16} \tilde{g}_i + \frac{10h}{16} \tilde{g}_{i+1}  - \frac{h}{16} \tilde{g}_{i+2})+R_z^i\\
     &:=hE_i[Z_{i+1}] + h \sum_{j=0}^{Kz} \gamma_{K_z,j}^{1} E_i[f(t_{i+j},Y_{i+j},Z_{i+j})\Delta W_{i+j}]-h \sum_{j=0}^{Kz} \gamma_{K_z,j}^{1} E_i[Z_{i+j}]+R_z^i,
 \end{align}
where $\gamma_{K_z,0}^{1}=\frac{7}{16}, \gamma_{K_z,1}^{1}=\frac{5}{8}, \gamma_{K_z,2}^{1}=-\frac{1}{16}.$\\

Using the second order polynomials we obtain
\begin{equation}
\begin{split}
 p_i(s)&=g_i(s-t_i) -\left(\frac{3}{2}g_i-2g_{i+1}+ \frac{1}{2}g_{i+2}\right) (s-t_i) /h \\
       &+ \left(\frac{1}{2}g_i-g_{i+1}+ \frac{1}{2}g_{i+2}\right) (s-t_i)^2 /h^2
\end{split}
\end{equation}
\begin{equation}
\begin{split}
 \tilde{p}_i(s)&=\tilde{g}_i(s-t_i) -\left(\frac{3}{2}\tilde{g}_i-2\tilde{g}_{i+1}+ \frac{1}{2}\tilde{g}_{i+2}\right) (s-t_i) /h \\
       &+ \left(\frac{1}{2}\tilde{g}_i-\tilde{g}_{i+1}+ \frac{1}{2}\tilde{g}_{i+2}\right) (s-t_i)^2 /h^2
\end{split}
\end{equation}
whose integrations are given by
\begin{equation}
 \int_{t_i}^{t_{i+1}}p_i(s)ds= h\frac{5g_i+8g_{i+1}-g_{i+2}}{12}.
\end{equation}
\begin{equation}
 \int_{t_i}^{t_{i+1}}\tilde{p}_i(s)ds= h\frac{5\tilde{g}_i+8\tilde{g}_{i+1}-\tilde{g}_{i+2}}{12}.
\end{equation}
By using the second order polynomial we rewrite \eqref{eq:reference_y_01} as 
\begin{align}
  &0 = hE_i[Z_{i+1}] + \frac{5h}{12} g_i + \frac{2h}{3} g_{i+1}  - \frac{h}{12} g_{i+2} - (\frac{5h}{12} \tilde{g}_i + \frac{2h}{3} \tilde{g}_{i+1}  - \frac{h}{12} \tilde{g}_{i+2})+R_z^i\\
     &:=hE_i[Z_{i+1}] + h \sum_{j=0}^{Kz} \gamma_{K_z,j}^{1} E_i[f(t_{i+j},Y_{i+j},Z_{i+j})\Delta W_{i+j}]-h \sum_{j=0}^{Kz} \gamma_{K_z,j}^{1} E_i[Z_{i+j}]+R_z^i,
 \end{align}
where $\gamma_{K_z,0}^{1}=\frac{5}{12}, \gamma_{K_z,1}^{1}=\frac{2}{3}, \gamma_{K_z,2}^{1}=-\frac{1}{12}.$

\item $K_z=3:$ for $K_z \geq 3$ we will use the Not-a-knot cubic spline.\\

For $\ts^{z_1}_{t_i,0}(s), \ts^{z_2}_{t_i,0}(s), s\in[t_i,t_{i+1}]: $
\begin{align*}
 a^{z_1}_0 &= g_i, b^{z_1}_0 =-(11g_i-18g_{i+1} + 9g_{i+2} - 2g_{i+3})/6h\\
 c^{z_1}_0 &=(2g_i-5g_{i+1} + 4g_{i+2} - g_{i+3})/2h^2, d^{z_1}_0 =-(g_i-3g_{i+1} + 3g_{i+2}-g_{i+3})/6h^3
\end{align*}
\begin{align*}
 a^{z_2}_0 &= \tilde{g}_i, b^{z_2}_0 =-(11\tilde{g}_i-18\tilde{g}_{i+1} + 9\tilde{g}_{i+2} - 2\tilde{g}_{i+3})/6h\\
 c^{z_2}_0 &=(2\tilde{g}_i-5\tilde{g}_{i+1} + 4\tilde{g}_{i+2} - \tilde{g}_{i+3})/2h^2, d^{z_2}_0 =-(\tilde{g}_i-3\tilde{g}_{i+1} + 3\tilde{g}_{i+2}-g_{i+3})/6h^3
\end{align*}
Thus, \eqref{eq:disBSDEZ_06} can be rewritten as
\begin{align}
  &0 =h E_i[Z_{i+1}] + \frac{3h}{8} g_i + \frac{19h}{24} g_{i+1}  - \frac{5h}{24} g_{i+2} +  \frac{h}{24} g_{i+3} \nonumber\\
       & - (\frac{3h}{8} \tilde{g}_i + \frac{19h}{24} \tilde{g}_{i+1}  - \frac{5h}{24} \tilde{g}_{i+2} +\frac{h}{24} \tilde{g}_{i+3})+R_z^i\\
     &:=hE_i[Z_{i+1}] + h \sum_{j=0}^{Kz} \gamma_{K_z,j}^{1} E_i[f(t_{i+j},Y_{i+j},Z_{i+j})\Delta W_{i+j}]-h \sum_{j=0}^{Kz} \gamma_{K_z,j}^{1} E_i[Z_{i+j}]+R_z^i,\label{eq:ex_z}
 \end{align}
where $\gamma_{K_z,0}^{1}=\frac{3}{8}, \gamma_{K_z,1}^{1}=\frac{19}{24}, \gamma_{K_z,2}^{1}=-\frac{5}{24}, \gamma_{K_z,3}^{1}=\frac{1}{24}.$
\end{itemize}
The coefficients for $1\leq K_z \leq 6$ are reported in Table \ref{tab: coeffs_z}.
\begin{table}[h]
	\centering
    \begin{tabular}{|c|c|c|c|c|c|c|c|}
    \hline
    $K_z$ & \multicolumn{7}{|c|}{$\gamma_{K,j}^1$}  \\ \hline
          &$j=0$ & $j=1$ &$j=2$&$j=3$ &$j=4$ &$j=5$ &$j=6$\\ \hline
    1     & $\frac{1}{2}$  & $\frac{1}{2}$ & & & & & \\    \hline
    2 (Second Order Polynomial)    & $\frac{5}{12}$  & $\frac{2}{3}$ &$-\frac{1}{12}$ & & & & \\    \hline
    2 (Natural Cubic Spline )     & $\frac{7}{16}$  & $\frac{5}{8}$ &$-\frac{1}{16}$ & & & & \\    \hline
    3     & $\frac{3}{8}$  & $\frac{19}{24}$ &$-\frac{5}{24}$ &$\frac{1}{24}$ &  & &\\    \hline
    4     & $\frac{35}{96}$  & $\frac{5}{6}$ &$-\frac{13}{48}$ &$\frac{1}{12}$ &$-\frac{1}{96}$ & & \\    \hline
    5     & $\frac{131}{360}$  & $\frac{151}{180}$ &$-\frac{103}{360}$ &$\frac{37}{360}$ &$-\frac{1}{45}$ &$\frac{1}{360}$ & \\    \hline
    6     & $\frac{163}{448}$  & $\frac{47}{56}$ &$-\frac{129}{448}$ &$\frac{3}{28}$ &$-\frac{37}{1344}$ &$\frac{1}{168}$ &$-\frac{1}{1344}$ \\    \hline
    \end{tabular}
    \caption{The coefficients $[\gamma_{K_z,j}^{1}]_{j=0}^{K_z}$ for $K_z=1,2,\cdots,6.$}
    \label{tab: coeffs_z}
\end{table}
Note that $\Delta W_{t_i}=0$ and $E_i[Z_i]=Z_i,$ based on the calculations above we can obtain the reference equations of
the BSDEs as
 \begin{align}
  Y_i &= E_i[Y_{i+K_y}] + h K_y\sum_{j=0}^{Ky} \gamma_{K_y,j}^{K_y} E_i[f(t_{i+j},Y_{i+j},Z_{i+j})]+R_y^i,\label{eq:r_y}\\
  Z_i&=\left(E_i[Z_{i+1}] + \sum_{j=1}^{Kz} \gamma_{K_z,j}^{1} E_i[f(t_{i+j},Y_{i+j},Z_{i+j})\Delta W_{i+j}]- \sum_{j=1}^{Kz} \gamma_{K_z,j}^{1} E_i[Z_{i+j}]\right)/ \gamma_{K_z,0}^{1}+R_z^i,\label{eq:r_z}
   \end{align}
%\eqref{eq:disBSDEZ_05} is thus an explicit scheme.   
where $Y_i=\left(Y_i^{1},Y_i^{2},\cdots,Y_i^{m}\right)^{\top}$, $Z_i=\left(Z^{\tilde{m},\tilde{d}}_i\right)_{m\times d},$
$\Delta W_{i+j}=(W_{i+j}^1, W_{i+j}^2,\cdots,W_{i+j}^d)^{\top}-(W_i^1, W_i^2,\cdots,W_i^d)^{\top},$
$R_{y}^{i}=\left(R_{y}^{i,1},R_{y}^{i,2},\cdots,R_{y}^{i,{m}}\right)^{\top}$ and $R_{z}^{i}=\left(R_{z}^{i,{\tilde{m},\tilde{d}}}\right)_{m\times d}.$
It is easy to see that \eqref{eq:r_y} is implicit, and \eqref{eq:r_z} is always explicit for solving $Z_i.$
One can show that estimates for the local error terms $R^i_y$ and $R^i_z$ (componentwise in
\eqref{eq:residuals_y} and \eqref{eq:residuals_z}) are given by
\begin{equation}
 |R_y^i|=\mathcal{O}(h^{5}),\quad |R_z^i|=\mathcal{O}(h^{5})
\end{equation}
provided that the generator function $f$ and the terminal function $g$ are smooth.
It is worth noting that $R_z^i$ will be divided by $h$ for solving $Z_i,$ see e.g., \eqref{eq:ex_z},
one might set $K_z=K_y+1$ in order to balance the local truncation errors.
\section{A stable multistep discretization scheme}
In this Section we present a stable multistep scheme fully discrete in time and space.
\subsection{The Semi-discretization in time}
We denote $Y^i=\left(Y^{1,i},Y^{2,i},\cdots,Y^{m,i}\right)^{\top}$ and $Z^i=\left(Z^{\tilde{m},\tilde{d},i}\right)_{m\times d}$ as the approximations to $Y_i$ and $Z_i,$ namely
at the time $t_i$ in the reference equations, respectively.
Furthermore, we have $W_i=(W_i^1, W_i^2,\cdots,W_i^d)^{\top},$ whereas all Brownian motions are independent with each other.
Since $Z_i$ is needed for computing $Y_i$ in our scheme, we thus need to consider the larger step size between $K_y$ and $K_z.$ Therefore, we define the number of time steps as $K=\max\left\{K_y, K_z \right\}.$
Suppose that the random variables $Y^{N_T-j}$ and $Z^{N_T-j}$ are given for $j=0,1,\cdots,K-1,$ then $Y^{i}$ and $Z^{i}$ can be found for $i=N_T-K,\cdots,0$ by
 \begin{align}
  Y^i &= E_i[Y^{i+K_y}] + h K_y\sum_{j=0}^{Ky} \gamma_{K_y,j}^{K_y} E_i[f(t_{i+j},Y^{i+j},Z^{i+j})],\label{eq:r_y_01}\\
  Z^i&=\left(E_i[Z^{i+1}] + \sum_{j=1}^{Kz} \gamma_{K_z,j}^{1} E_i[f(t_{i+j},Y^{i+j},Z^{i+j})\Delta W^{\top}_{i+j}]- \sum_{j=1}^{Kz} \gamma_{K_z,j}^{1} E_i[Z^{i+j}]\right)/ \gamma_{K_z,0}^{1},\label{eq:r_z_01}
   \end{align}
We follow the methodologies used in \cite{Zhao2010} to check the stability.
We set the generator function $f=0$ and take the expectation $E[\cdot]$ on both sides of \eqref{eq:r_y_01}
\begin{equation}\label{eq:stable_y}
E[Y^i] = E[Y^{i+k}].
\end{equation}
Note that we have set $k=K_y$ in \eqref{eq:r_y_01}. We need to recall $k$ in \eqref{eq:stable_y} for a general stability analysis.
\eqref{eq:stable_y} indicates that reference equation of $Y_i$ is stable for any integers $1\leq k \leq K_y \leq N_T.$
Furthermore, in \eqref{eq:r_y_01} where $k=K_y,$ we have checked that $\sum_{j=0}^{K_y}\gamma_{K_y,j}^{K_y}=1$ for $1 \leq K_y \leq N_T.$

In a similar way to above, \eqref{eq:r_z} can be reformulated as
\begin{equation}\label{eq:stable_z}
 0=E[Z^{i+l}] - \sum_{j=1}^{Kz} \gamma_{K_z,j}^{l} E[Z^{i+j}],
\end{equation}
where $l$ is recalled substituting $1$ in \eqref{eq:r_z_01}. We see that \eqref{eq:stable_z} is a difference equation
of $Z^i,$ the characteristic polynomial of the backward difference equation \eqref{eq:stable_z} reads
\begin{equation}\label{eq:cf}
 p^l_{K_z}(\lambda)=\lambda^{K_z - l} - \sum_{j=1}^{Kz} \gamma_{K_z,j}^{l} \lambda^{K_z-j}.
\end{equation}
In order to have a stable reference equation of $Z^i,$ the roots of \eqref{eq:cf} must satisfy the following
condition:
\begin{itemize}
 \item The roots must be in the closed unit disc and the ones on the unit circle must be simple.
\end{itemize}
The values of $ \gamma_{K_z,j}^{1}$ have been given for $K_z=1,\cdots,6$ in Table \ref{tab: coeffs_z}.
In the same way as we obtained those values one can calculate the values of $ \gamma_{K_z,l}^{j}$ for $1<l\leq K_z\leq N_T$
and obtain the corresponding roots of \eqref{eq:cf}, see Table \ref{tab:roots}.
\begin{table}[h!]
	\centering
	\scalebox{0.55}{
 \begin{tabular}{|c|c|c|c|c|c|c|c|}
    \hline
    $K_z$ & $l$ & \multicolumn{6}{|c|}{Roots $\lambda^l_{K_z,j}$}  \\ \hline
          &     & $j=1$ &$j=2$&$j=3$ &$j=4$ &$j=5$ &$j=6$  \\ \hline
    1     &  1  & 1 & & & & &\\    \hline
    \multirow{2}{*}[-1em]{2} & 1 & 1 &\sr $-\frac{1}{5}\:\:\:(-\frac{1}{7}\text{ natural CS})$ & & & &\\ \cline{2-2}
    &2 &1 &\sr $-5 \:\:\:(-4.3333\text{ natural CS})$ & & & &\\ \hline
    \multirow{3}{*}[-1.6em]{3} & 1 & 1 &\sr $\frac{\sqrt{13}}{9} - \frac{2}{9}$ & $-\frac{\sqrt{13}}{9} - \frac{2}{9} $ & & &\\ \cline{2-2}
    & 2 & 1 &\sr $0$ & $-5$ & & &\\ \cline{2-2}
    & 3 & 1 &\sr $\sqrt{3}i - 2 $ & $-\sqrt{3}i - 2$ & & &\\ \hline
    \multirow{4}{*}[-2.3em]{4} & 1 & 1 &\sr $-0.82662$ & $0.14188 - 0.12014i $ &$0.14188 + 0.12014i$ & &\\ \cline{2-2}
    & 2 & 1 &\sr $0$ & $0$ &$-5$ & &\\ \cline{2-2}
    & 3 & 1 &\sr $-0.01244$ & $- 2.31196 + 1.40033i$ &$- 2.31196 - 1.40033i$ & &\\ \cline{2-2}
    & 4 & 1 &\sr $-3.93114$ & $- 0.53442 + 1.5851i$ &$- 0.53442 - 1.5851i$ & &\\ \hline
    \multirow{5}{*}[-2.3em]{5} & 1 & 1 &\sr $-0.89193$ & $0.20080$ &$0.06693 + 0.19529i$ &$0.06693 -	 0.19529i$ &\\ \cline{2-2}
    & 2 & 1 &\sr $0$ & $0$ &$0$ &$-5$ &\\ \cline{2-2}
    & 3 & 1 &\sr $-0.07259$ & $0.04667$ &$- 2.34069 - 1.31158i$ &$- 2.34069 + 1.31158i$ &\\ \cline{2-2}
    & 4 & 1 &\sr $-3.64370$ & $-0.00620$ &$- 0.57668 - 1.60195i$ &$- 0.57668 + 1.60195i$ &\\ \cline{2-2}
    & 5 & 1 &\sr $- 2.45215 + 0.06565i$ & $- 2.45215 - 0.06565i$ &$- 0.09849 - 1.50203i$ &$- 0.09849 + 1.50203i$ &\\ \hline
    \multirow{6}{*}[-2.3em]{6} & 1 & 1 &\sr $-0.91034$ & $- 0.01033 - 0.22612i $ &$- 0.01033 + 0.22612i $ &$0.18636 - 0.09543i $ &$0.18636 + 0.09543i $\\ \cline{2-2}
    & 2 & 1 &\sr $0$ & $0$ &$0$ &$0$ &$-5$\\ \cline{2-2}
    & 3 & 1 &\sr $-0.13432$ & $- 2.34031 + 1.29934i$ &$- 2.34031 - 1.29934i$ &$0.05126 + 0.06452i$ &$0.05126 - 0.06452i$\\ \cline{2-2}
    & 4 & 1 &\sr $-3.61188$ & $-0.04794$ &$0.03504$ &$- 0.58234 - 1.59752i$ &$- 0.58234 + 1.59752i$\\ \cline{2-2}
    & 5 & 1 &\sr $-3.00560$ & $-1.94659$ &$-0.00538$ &$0.09695 - 1.51077i$ &$0.09695 + 1.51077i$\\ \cline{2-2}
    & 6 & 1 &\sr $-3.38909$ & $-1.14732 + 1.07617i$ &$-1.14732 - 1.07617i$ &$0.44714 + 1.33772i$ &$0.44714 - 1.33772i$\\ \hline
    \end{tabular}}
 \caption{The roots of \eqref{eq:cf} for $K_z=1,2,\cdots,6$ and $l=1,\cdots,K_z$}
    \label{tab:roots}
\end{table}

Note that, for $K_y=1, 2, 3$ and $K_z=1, 2, 3,$ our reference equations (with second order polynomial for $K=2$) coincide with the reference equations proposed in
\cite{Zhao2010}, where Lagrange interpolating polynomials are employed. However, in \cite{Zhao2010}, the reference equation of $Y^i$ is stable only when
$K_y=1,2,3,4,5,6,7,9;$ and the reference equation of $Z^i$ is stable only when $K_z=1,2,3.$ As mentioned already, our both reference equations are generally stable,
namely for all $K_y\geq 1$ and $K_z\geq 1.$ This is to say that our method allows for considering more multi-time levels.
\subsection{Error analysis}
Due to the nested conditional expectations we still are confronted with a problem to perform error analysis for the proposed multi-step scheme.
In \cite{Zhao2010}, the authors have finished some error analysis for the multi-step semidiscrete scheme in one-dimensional case using the Lagrange interpolating
polynomials under several assumptions. In this section, we adopt their results to our multi-step scheme.
Throughout this section we assume that the functions $f$ and $g$ are bounded and smooth enough with bounded derivatives for a uniquely
existing solution. Furthermore, suppose that $f$ does not involve the variable $Z_t,$ i.e.,
\begin{equation}
 Y_t = \xi + \int_{t}^T f(s, Y_s)\,ds - \int_t^T Z_s\,dW_s,
\end{equation}
for which the reference equation read
 \begin{align}
  Y_i &= E_i[Y_{i+K_y}] + h K_y\sum_{j=0}^{Ky} \gamma_{K_y,j}^{K_y} E_i[f(t_{i+j},Y_{i+j})]+R_y^i,\label{eq:error_y}\\
  Z_i&=\left(E_i[Z_{i+1}] + \sum_{j=1}^{Kz} \gamma_{K_z,j}^{1} E_i[f(t_{i+j},Y_{i+j})\Delta W_{i+j}]- \sum_{j=1}^{Kz} \gamma_{K_z,j}^{1} E_i[Z_{i+j}]\right)/ \gamma_{K_z,0}^{1} + R_z^i/h,\label{eq:error_z}
   \end{align}
where the local truncation errors $R_y^i$ and $R_z^i$ are defined in \eqref{eq:residual_y} and \eqref{eq:residual_z}. 
And the corresponding multi-step scheme for $Y^i$ and $Z^i$ can be immediately written down from \eqref{eq:r_y_01} and \eqref{eq:r_z_01}. 
   \begin{lemma}
The local estimates of the local truncation errors in \eqref{eq:error_y} and \eqref{eq:error_z} satisfy 
 \begin{equation*}
 |R_y^i| \leq C h^{\min\{K_y+2, \,5\}}
 \quad
 |R_z^i| \leq C h^{\min\{K_z+2,\,5\}},
\end{equation*}
\iffalse
 \begin{equation*}
\left\{
\begin{array}{l}
 |R_y^i| \leq C h^{K_y+2},~\mbox{for}~K_y=1, 2, 3\\
 |R_y^i| \leq Ch^{5},~\mbox{for}~K_y > 3
 \end{array}
 \right.
 ,~
 \left\{
\begin{array}{l}
 |R_z^i| \leq C h^{K_z+2},~\mbox{for}~K_z=1, 2, 3\\
 |R_z^i| \leq Ch^{5},~\mbox{for}~K_z > 3
 \end{array}
 \right.
\end{equation*}
.\fi
where $C>0$ is a constant depending on $T, f, g$ and the derivatives of $f, g.$
\end{lemma}
The proof can be done directly by combining the proof of Lemma 3.2 in \cite{Zhao2009} and the fact that
not-a-knot cubic spline is fourth-order accurate.
\begin{theorem}\label{theorem:error_y}
 Suppose that the initial values satisfy
 %\begin{equation*}
 % \max_{N_T-K_y<i\leq N_T}E\left[\left|Y_{i} - Y^i\right|\right]=\mathcal{O}(h^{\min\{K_y+1,\,4\}})
 %\end{equation*}
 \begin{equation*}
\left\{
\begin{array}{l}
 \max_{N_T-K_y<i\leq N_T}E\left[\left|Y_{i} - Y^i\right|\right]=\mathcal{O}(h^{K_y+1}),~\mbox{for}~K_y=1, 2, 3\\
 \max_{N_T-K_y<i\leq N_T}E\left[\left|Y_{i} - Y^i\right|\right]=\mathcal{O}(h^4),~\mbox{for}~K_y > 3
 \end{array}
 \right.
\end{equation*}
 for sufficiently small time step $h$ it can be shown that
 \begin{equation}
\sup_{0\leq i \leq N_T} E\left[\left|Y_{i} - Y^i\right|\right] \leq Ch^{\min\{K_y+1,\, 4\}},
 \end{equation}
\iffalse
 \begin{equation*}
\left\{
\begin{array}{l}
 \sup_{0\leq i \leq N_T} E\left[\left|Y_{i} - Y^i\right|\right] \leq Ch^{K_y+1},~\mbox{for}~K_y=1, 2, 3\\
 \sup_{0\leq i \leq N_T} E\left[\left|Y_{i} - Y^i\right|\right] \leq Ch^{4},~\mbox{for}~K_y > 3
 \end{array}
 \right.
\end{equation*}
.\fi
where $C>0$ is a constant depending on $T, f, g$ and the derivatives of $f, g.$
\end{theorem}
The proof can be done directly by combining the proof of Theorem 1. in \cite{Zhao2010} and the fact that
not-a-knot cubic spline is fourth-order accurate.
\begin{theorem}\label{theorem:generalTheorem}
 Suppose that the initial values satisfy
 \begin{equation*}
\left\{
\begin{array}{l}
 \max_{N_T-K_z<i\leq N_T}E\left[\left|Z_{i} - Z^i\right|\right]=\mathcal{O}(h^{K_z}),~\mbox{for}~K_z=1, 2, 3\\
 \max_{N_T-K_z<i\leq N_T}E\left[\left|Z_{i} - Z^i\right|\right]=\mathcal{O}(h^3)~\mbox{for}~K_z > 3
 \end{array}
 \right.
\end{equation*}
and the condition on the initial values in Theorem \ref{theorem:error_y} is fulfilled.
For sufficiently small time step $h$ it can be shown that
\begin{equation*}
 \sup_{0\leq i \leq N_T} E\left[\left|Z_{i} - Z^i\right|\right] \leq Ch^{\min(K_y+1,\,K_z,\,3)},
\end{equation*}
\iffalse
\begin{equation*}
\left\{
\begin{array}{l}
 \sup_{0\leq i \leq N_T} E\left[\left|Z_{i} - Z^i\right|\right] \leq Ch^{\min(K_y+1,K_z)},~\mbox{for}~K_y=1, 2, K_z=1, 2, 3\\
 \sup_{0\leq i \leq N_T} E\left[\left|Z_{i} - Z^i\right|\right] \leq Ch^{3},~\mbox{for}~K_y>2, K_z > 3
 \end{array}
 \right.
\end{equation*}
.\fi
where $C>0$ is a constant depending on $T, f, g$ and the derivatives of $f, g.$
\end{theorem}
The proof can be done directly by combining the proof of Theorem 2. in \cite{Zhao2010} and the fact that
not-a-knot cubic spline is fourth-order accurate.

\subsection{The fully discretized scheme }
We have checked that \eqref{eq:r_y_01} and \eqref{eq:r_z_01} are stable in the time direction. To solve
$(Y^i,Z^i)$ numerically, next we consider the space discretization. We define firstly the partion of the one-dimensional $(\tilde{d}=d=1)$ real axis as
\begin{equation}
\mathcal{R}^{\tilde{d}}  =\left\{x_{\gamma}^{\tilde{d}}| x_{\gamma}^{\tilde{d}} \in \mathbb{R}, \gamma \in \mathbb{Z}, 
x_{\gamma}^{\tilde{d}}<x_{\gamma+1}^{\tilde{d}}, \lim_{i \to +\infty} x_{\gamma}^{\tilde{d}}=+\infty,\lim_{i \to -\infty} x_{\gamma}^{\tilde{d}}=-\infty \right\}.
\end{equation}
Thus, the partition of $d$-dimensional space $\mathcal{R}^d$ reads
\begin{equation}
 \mathcal{R}^{\tilde{d}}=\mathcal{R}^{1} \times \cdots \times \mathcal{R}^{\tilde{d}}\times \cdots \times \mathcal{R}^{d},
\end{equation}
where $\tilde{d}=1,2,\cdots,d.$ For simplicity of notation we will use $x_{\Gamma}=(x^1_{\gamma_1},x^2_{\gamma_2},\cdots,x^d_{\gamma_d})^{\top}$
for $\Gamma=(\gamma_1,\gamma_2,\cdots,\gamma_d)\in \mathbb{Z}^d.$
We use $y^{N_T-\lambda}_{\Gamma}$ and $z^{N_T-\lambda}_{\Gamma}$ to denote the values of random variables $Y^{N_T-\lambda}$ and $Z^{N_T-\lambda}$
at the points $x_{\Gamma}.$ Given these values for $\lambda=0,1,\cdots,K-1,$ we need to find $(y^{i}_{\Gamma},z^{i}_{\Gamma}), i=N_T-K,\cdots,0$ such that
\begin{align}
  y^i_{\Gamma} &= E^{x_{\Gamma}}_i[\hat{Y}^{i+K_y}] + h K_y\sum_{j=0}^{Ky} \gamma_{K_y,j}^{K_y} E^{x_{\Gamma}}_i[f(t_{i+j},\hat{Y}^{i+j},\hat{Z}^{i+j})],\label{eq:r_y_02}\\
  z^i_{\Gamma}&=\left(E^{x_{\Gamma}}_i[\hat{Z}^{i+1}] + \sum_{j=1}^{Kz} \gamma_{K_z,j}^{1} E^{x_{\Gamma}}_i[f(t_{i+j},\hat{Y}^{i+j},\hat{Z}^{i+j})\Delta W^{\top}_{i+j}]
  - \sum_{j=1}^{Kz} \gamma_{K_z,j}^{1} E^{x_{\Gamma}}_i[\hat{Z}^{i+j}]\right)/ \gamma_{K_z,0}^{1},\label{eq:r_z_02}
   \end{align}
where $E^{x_{\Gamma}}_i[\cdot]$ denotes the conditional expectation
under the $\sigma$-field $\mathcal{F}_t^{x_{\Gamma}}$ generated by $\{W_i=x_\Gamma\}.$
Correspondingly, $\hat{Y}^{i+j}$ and $\hat{Z}^{i+j}$ denote the functions of increment of Brownian
motion $Y^{i+j}(\Delta W_i)$ and $Z^{i+j}(\Delta W_i)$ with the fixed $\{W_i=x_\Gamma\}.$

To approximate the conditional expectations in \eqref{eq:r_y_02} and \eqref{eq:r_z_02} we employ the Gauss-Hermite quadrature formula.
For example, we compute $E^{x_{\Gamma}}_i[\hat{Y}^{i+K_y}]$ as
\begin{align}
 E^{x_{\Gamma}}_i[\hat{Y}^{i+K_y}]&=\frac{1}{(2K_y\pi h )^{d/2}}\int_{\mathbb{R}^d}\hat{Y}^{i+K_y}(s) \exp\left(-\frac{(s-x)^{\top}(s-x)}{2K_y h}\right)\,ds\\
                                  &\approx \frac{1}{(2K_y\pi h )^{d/2}}\int_{\mathbb{R}^d}\hat{y}^{i+K_y}(s) \exp\left(-\frac{(s-x)^{\top}(s-x)}{2K_y h}\right)\,ds\\
                                   &\approx \frac{1}{\pi^{\frac{d}{2}}}\sum_{\Lambda=1}^{L}\omega_{\Lambda}\hat{y}^{i+K_y}(x_{\Gamma}+\sqrt{2K_y h}a_{\Lambda})\label{eq:lroot}\\
                                   &:= \hat{E}^{x_{\Gamma}}_i[\hat{Y}^{i+K_y}],
\end{align}
where $\hat{y}^{i+K_y}(s)$ are interpolating values at the space points $s$ based on $y_{\Gamma}^{i+K_y}$ at a finite number
of the space grid points $x_{\Gamma}$ near $s,\, \Lambda=(\lambda_1, \lambda_2, \cdots, \lambda_d),\, 
\omega_{\Lambda}=\prod_{\tilde{d}=1}^{d}\omega_{\lambda_{\tilde{d}}},\, a_{\Lambda}=(a_{\lambda_1}, a_{\lambda_2}, \cdots, a_{\lambda_d}),\,
\sum_{\Lambda=1}^{L}=\sum_{\lambda_1=1,\cdots,\lambda_d=1}^{L, \cdots, L}.$ For the weights $\omega_{\Lambda}$ and the roots $a_{\Lambda}$
we refer to e.g., \cite{Abramowitz1972}. The approximations of the other conditional expectations in \eqref{eq:r_y_02} and \eqref{eq:r_z_02} can be
done similarly. Finally, by considering these approximations we rewrite \eqref{eq:r_y_02} and \eqref{eq:r_z_02} as
\begin{align}
  y^i_{\Gamma} &= \hat{E}^{x_{\Gamma}}_i[\hat{Y}^{i+K_y}] + h K_y\sum_{j=0}^{Ky} \gamma_{K_y,j}^{K_y} \hat{E}^{x_{\Gamma}}_i[f(t_{i+j},\hat{Y}^{i+j},\hat{Z}^{i+j})],\label{eq:r_y_03}\\
  z^i_{\Gamma}&=\left(\hat{E}^{x_{\Gamma}}_i[\hat{Z}^{i+1}] + \sum_{j=1}^{Kz} \gamma_{K_z,j}^{1} \hat{E}^{x_{\Gamma}}_i[f(t_{i+j},\hat{Y}^{i+j},\hat{Z}^{i+j})\Delta W^{\top}_{i+j}]
  - \sum_{j=1}^{Kz} \gamma_{K_z,j}^{1} \hat{E}^{x_{\Gamma}}_i[\hat{Z}^{i+j}]\right)/ \gamma_{K_z,0}^{1}.\label{eq:r_z_03}
   \end{align}
We observe that the computations at each space grid point are independent, which can be thus parallelized.
Usually, only the values of $y^{N_T}_{\Gamma}$ and $z^{N_T}_{\Gamma}$ are known because of the terminal condition.
However, for a $K$-step scheme we need to know the support values of $y^{N_T-j}_{\Gamma}$ and $z^{N_T-j}_{\Gamma}, j=0,\cdots,K-1.$
One can use the following two ways to deal with this problem:
before running the multi-step scheme, we choose a quite smaller $h$ and run one-step scheme until $N_T-K;$
Alternatively, one can prepare these initial values ``iteratively'', namely we compute $y^{N_T-1}_{\Gamma}$ and $z^{N_T-1}_{\Gamma}$
based on $y^{N_T}_{\Gamma}$ and $z^{N_T}_{\Gamma}$ with $K=1,$ and the compute $y^{N_T-2}_{\Gamma}$ and $z^{N_T-2}_{\Gamma}$
based on  $y^{N_T}_{\Gamma}, y^{N_T-1}_{\Gamma}, z^{N_T}_{\Gamma}, z^{N_T-1}_{\Gamma}$ with $K=2$ and so on.
Notice that we are faced with a computational complexity problem for solving high-dimensional problem, since the number of the Gauss-Hermite
quadrature points grows exponentially with the dimension $d.$

\section{Numerical experiments}
In this section we use some numerical examples to show the high effectiveness and accuracy of our scheme for solving the BSDEs.
We choose the truncated domain for the Brownian motion to be $[-8, 8]^d,$ and the degree of the Hermite polynomial (see $L$ in \eqref{eq:lroot} )to be $8.$
Note that, for $L=8,$ the quadrature error is so small that it cannot affect the convergence rate.
We use the Newton-Raphson method to implicitly solve \eqref{eq:r_y_03}. For the interpolation method we apply cubic spline interpolation which is a
fourth-order accurate, namely $(\Delta x)^4.$ In order to be able to estimate the convergence rate in time, we adjust the space step size $\Delta x$ according
to the time step size $h$ such that $(\Delta x)^4=(h )^{q+1}$ with $q=\min\{K_y+1, K_z\}.$
In the general case (the generator $f$ depends on both $Y_t$ and $Z_t$),
from Theorem \ref{theorem:generalTheorem} we know that $q$ is only limited to 3, since not-a-knot cubic spline is maximal fourth-rate accurate. This is to say that we always take $q=3$ when $\min\{K_y+1, K_z\}\geq 3.$
However, when the generator $f$ does not involve the component $Z_t,$ the approximation for $Y_t$ can reach fourth-order accurate, see Theorem \ref{theorem:error_y}.
For this case, $q$ is allowed to be $4$ when $\min\{K_y+1, K_z\}\geq 4.$

Generally, only $Y_{N_T}$ and $Z_{N_T}$ are known analytically. However, as mentioned before, for a $K$-step scheme we need to know
$y^{N_T-j}_{\Gamma}$ and $z^{N_T-j}_{\Gamma}, j=1,\cdots,K-1$ as initial values as well. To obtain these initial values, we start with $K=1$
and choose a extremely small time step size $h.$ Because the largest number of steps in our experiments is $K=6,$ we start thus with $N_T=8.$
In our computation we have used parallel computing using Python's multiprocessing module.
Note that a GPU-based parallelism will be much more cost-effective, which is left as a future work.

As mentioned before, our algorithm coincides with the algorithm proposed in \cite{Zhao2010} for $K_y=1, 2, 3$ and $K_z=1, 2, 3.$
In \cite{Zhao2010}, the authors have compared the multi-step scheme to the implicit Euler scheme \cite{Zhao2009} and the $\theta$-scheme \cite{Zhao2006}.
For these implicit Euler scheme and $\theta$-scheme, they have considered both the Monte-Carlo method and the Gaussian quadrature for approximating the conditional expectations.
Therefore, we will not do any comparison with other methods, for this we refer \cite{Zhao2010}.
In our numerical examples we will demonstrate higher effectiveness and accuracy of our scheme, which allows for more than $3$-step scheme, namely $K>3.$

\paragraph{Example 1}
The first example reads
\begin{equation*}
\begin{cases}
&-dY_t=-\frac{5}{8}Y_t \,dt-Z_t\,dW_t,\\
&Y_T=\exp(W_T/2+T/2),
\end{cases}
\end{equation*}
with the analytic solution
\begin{equation*}
\begin{cases}
&Y_t=\exp(W_t/2+t/2),\\
&Z_t=\exp(W_t/2+t/2)/2.
\end{cases}
\end{equation*}
The exact solution of $(Y_0,Z_0)$ is thus $\left(1,\frac{1}{2}\right).$ Obviously, in this example, the generator $f$ does not depend on
$Z_t.$ We thus choose $q=\min\{K_y+1, K_z\} < 4$ and keep $q=4$ when $\min\{K_y+1, K_z\}\geq 4.$ This is to say that the value of $q$ is exactly
the theoretical convergence order for the $Y$-component solver. For the $Z$-component, the theoretical convergence order of our scheme is $\min\{K_y+1, K_z\}$
but limited by $3$ due to Theorem \ref{theorem:generalTheorem}. The corresponding numerical results and estimated convergence rates are reported in Table \ref{tab:exp01y}
and \ref{tab:exp01z}. For $K=1, \cdots, 4,$ we have considered many combinations with the different values of $K_y, K_z$ and the corresponding values of $q.$
The results of these combinations are also similar for $K\geq 5.$ Therefore, for $K=5, 6$ we only report the results for $K_y=K_z=5, 6$ which are sufficient to
show the benefit from a higher number of multi-step.
\begin{table}[h!]
\centering
\begin{tabular}{|c|c|c|c|c|c|c|}
\hline
\multicolumn{7}{|c|}{$|Y_0-y_0^0|$}  \\ \hline
 & $N_T=8$ & $N_T=16$ & $N_T=32$ & $N_T=64$ & $N_T=128$ & CR \\
\hline
$K_y=1, K_z=1, q=1$ & 3.40e-04 & 8.90e-05 & 2.48e-05 & 7.37e-06 & 2.48e-06 & 1.78 \\
\hline
$K_y=1, K_z=2, q=2$ & 3.19e-04 & 7.96e-05 & 2.00e-05 & 5.00e-06 & 1.25e-06 & 2.00 \\
\hline
$K_y=2, K_z=1, q=1$ & 6.26e-06 & 2.81e-06 & 1.46e-06 & 7.16e-07 & 3.69e-07 & 1.01 \\
\hline
$K_y=2, K_z=2, q=2$ & 8.79e-07 & 3.24e-07 & 4.57e-08 & 8.83e-09 & 4.28e-09 & 2.06 \\
\hline
$K_y=2, K_z=3, q=3$ & 2.05e-07 & 1.16e-08 & 2.03e-09 & 1.38e-10 & 3.11e-11 & 3.18 \\
\hline
$K_y=3, K_z=1, q=1$ & 7.33e-07 & 2.06e-07 & 1.75e-07 & 8.88e-08 & 5.67e-08 & 0.86 \\
\hline
$K_y=3, K_z=2, q=2$ & 6.60e-07 & 8.09e-08 & 2.55e-08 & 8.35e-09 & 1.33e-09 & 2.12 \\
\hline
$K_y=3, K_z=3, q=3$ & 2.30e-07 & 2.52e-08 & 1.79e-09 & 2.58e-10 & 2.29e-11 & 3.32 \\
\hline
$K_y=3, K_z=4, q=4$ & 1.99e-07 & 1.77e-08 & 1.07e-09 & 7.05e-11 & 4.50e-12 & 3.88 \\
\hline
$K_y=4, K_z=1, q=1$ & 3.23e-07 & 5.36e-07 & 2.54e-07 & 1.42e-07 & 6.64e-08 & 0.64\\
\hline
$K_y=4, K_z=2, q=2$ & 5.11e-07 & 8.37e-08 & 3.77e-08 & 1.55e-09 & 1.68e-09 & 2.23 \\
\hline
$K_y=4, K_z=3, q=3$ & 1.54e-07 & 1.50e-08 & 9.64e-10 & 1.49e-10 & 8.19e-12 & 3.51 \\
\hline
$K_y=4, K_z=4, q=4$ & 1.54e-07 & 9.29e-09 & 5.59e-10 & 3.40e-11 & 2.04e-12 & 4.05 \\
\hline
$K_y=4, K_z=5, q=4$ & 1.54e-07 & 9.29e-09 & 5.59e-10 & 3.40e-11 & 2.04e-12 & 4.05 \\
\hline
$K_y=5, K_z=5, q=4$ & 6.48e-08 & 7.06e-09 & 4.12e-10 & 2.54e-11 & 1.66e-12 & 3.86 \\
\hline
$K_y=6, K_z=6, q=4$ & 6.60e-08 & 3.81e-09 & 3.21e-10 & 1.92e-11 & 1.32e-12 & 3.89 \\
\hline
\end{tabular}
\caption{Errors and convergence rates for Example 1, $T=1$}\label{tab:exp01y}
\end{table}
\begin{table}[h!]
\centering
\begin{tabular}{|c|c|c|c|c|c|c|}
\hline
\multicolumn{7}{|c|}{$|Z_0-z_0^0|$}  \\ \hline
 & $N_T=8$ & $N_T=16$ & $N_T=32$ & $N_T=64$ & $N_T=128$ & CR \\
\hline
$K_y=1, K_z=1, q=1$ & 1.71e-02 & 8.52e-03 & 4.25e-03 & 2.12e-03 & 1.06e-03 & 1.00 \\
\hline
$K_y=1, K_z=2, q=2$ & 8.50e-04 & 2.24e-04 & 5.76e-05 & 1.46e-05 & 3.67e-06 & 1.97 \\
\hline
$K_y=2, K_z=1, q=1$ & 1.72e-02 & 8.54e-03 & 4.26e-03 & 2.12e-03 & 1.06e-03 & 1.00 \\
\hline
$K_y=2, K_z=2, q=2$ & 7.89e-04 & 2.09e-04 & 5.37e-05 & 1.36e-05 & 3.42e-06 & 1.96 \\
\hline
$K_y=2, K_z=3, q=3$ & 4.17e-05 & 6.02e-06 & 8.03e-07 & 1.04e-07 & 1.32e-08 & 2.91 \\
\hline
$K_y=3, K_z=1, q=1$ & 1.72e-02 & 8.54e-03 & 4.26e-03 & 2.12e-03 & 1.06e-03 & 1.00 \\
\hline
$K_y=3, K_z=2, q=2$ & 7.89e-04 & 2.09e-04 & 5.37e-05 & 1.36e-05 & 3.42e-06 & 1.96 \\
\hline
$K_y=3, K_z=3, q=3$ & 4.16e-05 & 6.02e-06 & 8.03e-07 & 1.04e-07 & 1.32e-08 & 2.91 \\
\hline
$K_y=3, K_z=4, q=4$ & 1.98e-05 & 3.24e-06 & 4.59e-07 & 6.10e-08 & 7.84e-09 & 2.83 \\
\hline
$K_y=4, K_z=1, q=1$ & 1.72e-02 & 8.54e-03 & 4.26e-03 & 2.12e-03 & 1.06e-03 & 1.00 \\
\hline
$K_y=4, K_z=2, q=2$ & 7.89e-04 & 2.09e-04 & 5.37e-05 & 1.36e-05 & 3.42e-06 & 1.96 \\
\hline
$K_y=4, K_z=3, q=3$ & 4.17e-05 & 6.02e-06 & 8.03e-07 & 1.04e-07 & 1.32e-08 & 2.91 \\
\hline
$K_y=4, K_z=4, q=4$ & 1.98e-05 & 3.25e-06 & 4.60e-07 & 6.10e-08 & 7.90e-09 & 2.83 \\
\hline
$K_y=4, K_z=5, q=4$ & 1.67e-05 & 3.34e-06 & 5.00e-07 & 6.77e-08 & 1.30e-08 & 2.83 \\
\hline
$K_y=5, K_z=5, q=4$ & 1.67e-05 & 3.34e-06 & 4.99e-07 & 6.77e-08 & 1.10e-08 & 2.68 \\
\hline
$K_y=6, K_z=6, q=4$ & 1.29e-05 & 2.93e-06 & 4.61e-07 & 6.39e-08 & 1.60e-10 & 3.81 \\
\hline
\end{tabular}
\caption{Errors and convergence rates for Example 1, $T=1$}\label{tab:exp01z}
\end{table}

\begin{figure}[htbp!]
 \centering
 \begin{subfigure}[b]{0.47\textwidth}
 \includegraphics[width=\textwidth]{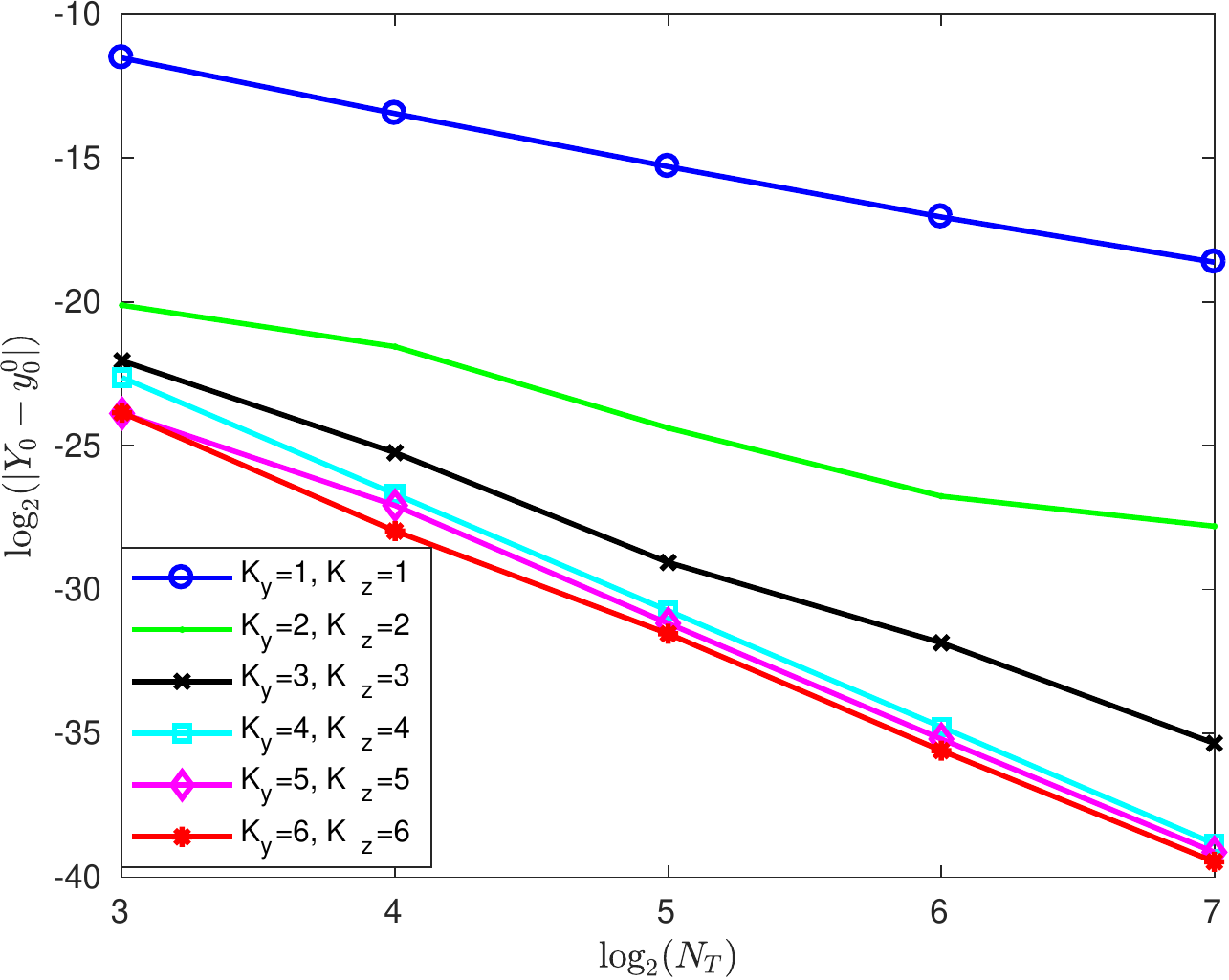}
 \subcaption{$Y$-component}
\end{subfigure}
 ~ 
 \begin{subfigure}[b]{0.47\textwidth}
 \includegraphics[width=\textwidth]{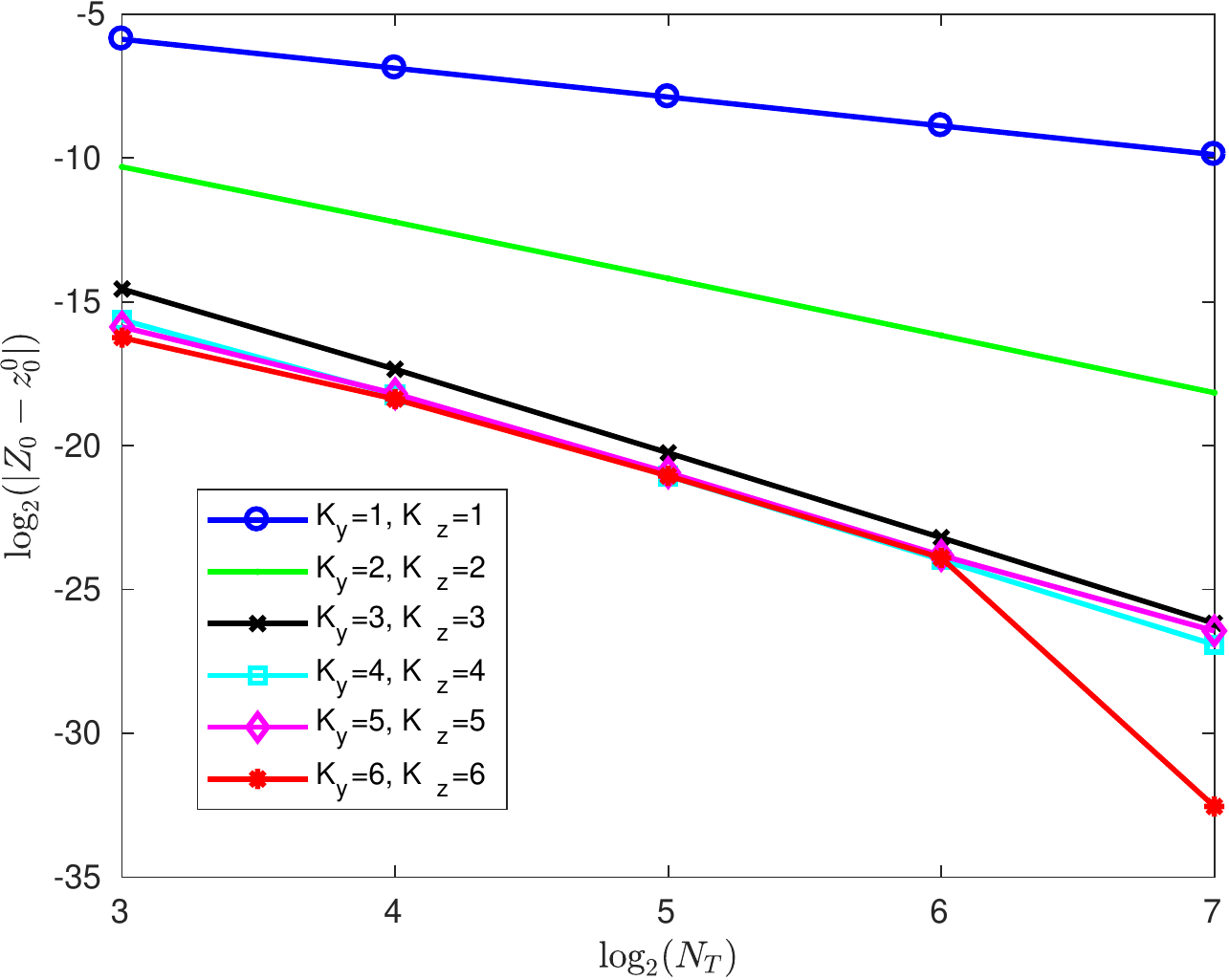}
 \subcaption{$Z$-component}
\end{subfigure}
 \caption{Plots of $\log_2\left(|Y_0-y_0^0|\right)$ and $\log_2\left(|Z_0-z_0^0|\right)$ with respect to $\log_2(N_T)$ for $K=1, \cdots 6$ for Example 1.}\label{fig:example_01}
\end{figure}

By a columnwise comparison we see that the approximation errors reduce mostly with the increasing number of steps, $K_y$ and $K_z.$
We have obtained $10^{-8}$ for approximating $Y_t$ already with $N_T=8,$ namely $h=\frac{1}{8}.$ The estimated convergence rates\footnote{Estimated by using linear squares fitting.} (CR) for both of
$Y_t$ and $Z_t$ are consistent with the theoretical results explained before, if we ignore the quadrature and interpolation errors which can cause
a slightly smaller estimated convergence rate. 
In Table \ref{tab:exp01z} we even observe a better CR than the theoretical result for $K_y=K_z=6.$ We display the plots of $\log_2\left(|Y_0-y_0^0|\right)$ and $\log_2\left(|Z_0-z_0^0|\right)$
with respect to $\log_2(N_T)$ in Figure \ref{fig:example_01}.

For this example, we also run our algorithm separately (without computing the $Z$-component) for solving the $Y$-component with smaller space step size $\Delta x$ (higher value of $q$). 
For $K_y\geq 4,$ we compare the numerical solutions computed with $q=4, \cdots, K_y+1.$
\begin{table}[h!]
\centering
\begin{tabular}{|c|c|c|c|c|c|c|}
\hline
\multicolumn{7}{|c|}{$|Y_0-y_0^0|$}  \\ \hline
 & $N_T=8$ & $N_T=16$ & $N_T=32$ & $N_T=64$ & $N_T=128$ & CR \\
\hline
$K_y=4, q=4$ & 1.54e-07 & 9.29e-09 & 5.59e-10 & 3.40e-11 & 2.04e-12 & 4.05 \\
\hline
$K_y=4, q=5$ & 1.53e-07 & 8.85e-09 & 5.30e-10 & 3.23e-11 & 2.00e-12 & 4.05 \\
\hline
\hline
$K_y=5, q=4$ & 6.48e-08 & 7.06e-09 & 4.12e-10 & 2.54e-11 & 1.66e-12 & 3.86 \\
\hline
$K_y=5, q=5$ & 6.24e-08 & 6.73e-09 & 4.03e-10 & 2.44e-11 & 1.63e-12 & 3.86 \\
\hline
$K_y=5, q=6$ & 6.21e-08 & 6.71e-09 & 4.02e-10 & 2.44e-11 & 1.63e-12 & 3.86 \\
\hline
\hline
$K_y=6, q=4$ & 6.60e-08 & 3.81e-09 & 3.21e-10 & 1.92e-11 & 1.32e-12 & 3.89\\
\hline
$K_y=6, q=5$ & 6.53e-08 & 3.62e-09 & 3.10e-10 & 1.87e-11 & 1.25e-12 & 3.89 \\
\hline
$K_y=6, q=6$ & 6.50e-08 & 3.62e-09 & 3.09e-10 & 1.87e-11 & 1.25e-12 & 3.89 \\
\hline
$K_y=6, q=7$ & 6.49e-08 & 3.62e-09 & 3.09e-10 & 1.87e-11 & 1.25e-12 & 3.89 \\
\hline
\end{tabular}
\caption{Errors and convergence rates for Example 1, where $y_0^0$ is separately computed for different higher values of $q$ and $T=1.$}\label{tab:exp01q}
\end{table}
The reported results in Table \ref{tab:exp01q} have shown clearly that there is almost no benefit to setting $q=K_y+1$ when $K_y+1>4,$ i.e.,
we only need to keep $q=4$ for $K_y+1 > 4.$ We emphasise again that the generator $f$ does not depends on $Z$-component in this example.
In general, this experiment clarifies that we should set $q=\min\{K_y+1, K_z\}<4$ and keep $q=3$ for $\min\{K_y+1, K_z\}\geq 4,$ the value of $q$
is thus the theoretical convergence order, see Theorem \ref{theorem:generalTheorem}.
\paragraph{Example 2}
For the second example we consider the nonlinear BSDE (taken from \cite{Zhao2010}) 
 \begin{equation*}
\begin{cases}
&-dY_t=\frac{1}{2}[\exp(t^2)-4tY_t-3\exp(t^2-Y_t\exp(-t^2))+Z_t^2\exp(-t^2)]\,dt-Z_t\,dW_t,\\
&Y_T=\ln(\sin W_T+3)\exp(T^2),
\end{cases}
\end{equation*}
with the analytic solution
\begin{equation*}
\begin{cases}
&Y_t=\ln\left(\sin{W_t}+3\right)\exp(t^2),\\
&Z_t=\exp(t^2)\frac{\cos{W_t}}{\sin{W_t}+3}.
\end{cases}
\end{equation*}
The exact solution of $(Y_0,Z_0)$ is then $\left(\ln(3),\frac{1}{3}\right).$
In this example, the generator $f$ is nonlinear and depends on $t, Y_t$ and $Z_t.$
Thus, from Theorem \ref{theorem:generalTheorem} we see that the theoretical
convergence order of our scheme for solving both $Y$ and $Z$ is $\min\{K_y+1, K_z\}$ but limited by $3.$
As clarified before, the used values of $q$ in both Table \ref{tab:exp02y}, \ref{tab:exp02z} are the values of
corresponding theoretical convergence order.
\begin{table}[h!]
\centering
\begin{tabular}{|c|c|c|c|c|c|c|}
\hline
\multicolumn{7}{|c|}{$|Y_0-y_0^0|$}  \\ \hline
  & $N_T=8$ & $N_T=16$ & $N_T=32$ & $N_T=64$ & $N_T=128$ & CR \\
\hline
$K_y=1, K_z=1, q=1$ & 2.72e-02 & 9.69e-03 & 3.87e-03 & 1.70e-03 & 7.87e-04 & 1.27 \\
\hline
$K_y=1, K_z=2, q=2$ & 1.40e-02 & 3.41e-03 & 8.43e-04 & 2.10e-04 & 5.22e-05 & 2.02 \\
\hline
$K_y=2, K_z=1, q=1$ & 1.17e-02 & 5.79e-03 & 2.89e-03 & 1.45e-03 & 7.24e-04 & 1.00 \\
\hline
$K_y=2, K_z=2, q=2$ & 1.38e-03 & 4.60e-04 & 1.27e-04 & 3.33e-05 & 8.47e-06 & 1.85 \\
\hline
$K_y=2, K_z=3, q=3$ & 6.39e-04 & 8.51e-05 & 1.13e-05 & 1.48e-06 & 1.89e-07 & 2.93 \\
\hline
$K_y=3, K_z=1, q=1$ & 1.05e-02 & 5.76e-03 & 2.87e-03 & 1.44e-03 & 7.22e-04 & 0.97 \\
\hline
$K_y=3, K_z=2, q=2$ & 1.44e-03 & 4.55e-04 & 1.27e-04 & 3.32e-05 & 8.48e-06 & 1.86 \\
\hline
$K_y=3, K_z=3, q=3$ & 5.34e-04 & 9.44e-05 & 1.19e-05 & 1.53e-06 & 1.92e-07 & 2.88 \\
\hline
$K_y=3, K_z=4, q=3$ & 2.33e-04 & 5.17e-05 & 6.55e-06 & 8.89e-07 & 1.13e-07 & 2.79 \\
\hline
$K_y=4, K_z=1, q=1$ & 1.19e-02 & 5.82e-03 & 2.89e-03 & 1.45e-03 & 7.23e-04 & 1.01 \\
\hline
$K_y=4, K_z=2, q=2$ & 1.38e-03 & 4.63e-04 & 1.28e-04 & 3.33e-05 & 8.48e-06 & 1.85 \\
\hline
$K_y=4, K_z=3, q=3$ & 6.60e-04 & 8.63e-05 & 1.14e-05 & 1.48e-06 & 1.89e-07 & 2.94 \\
\hline
$K_y=4, K_z=4, q=3$ & 3.49e-04 & 4.29e-05 & 6.04e-06 & 8.31e-07 & 1.10e-07 & 2.90 \\
\hline
$K_y=4, K_z=5, q=3$ & 3.33e-04 & 4.14e-05 & 6.18e-06 & 8.90e-07 & 1.21e-07 & 2.84 \\
\hline
$K_y=5, K_z=5, q=3$ & 1.13e-04 & 3.59e-05 & 5.81e-06 & 8.67e-07 & 1.20e-07 & 2.51 \\
\hline
$K_y=6, K_z=6, q=3$ & 8.55e-05 & 2.13e-05 & 4.75e-06 & 7.70e-07 & 1.11e-07 & 2.40 \\
\hline
\end{tabular}
\caption{Errors and convergence rates for Example 2, $T=1$}\label{tab:exp02y}
\end{table}

\begin{figure}[htbp!]
 \centering
 \begin{subfigure}[b]{0.47\textwidth}
 \includegraphics[width=\textwidth]{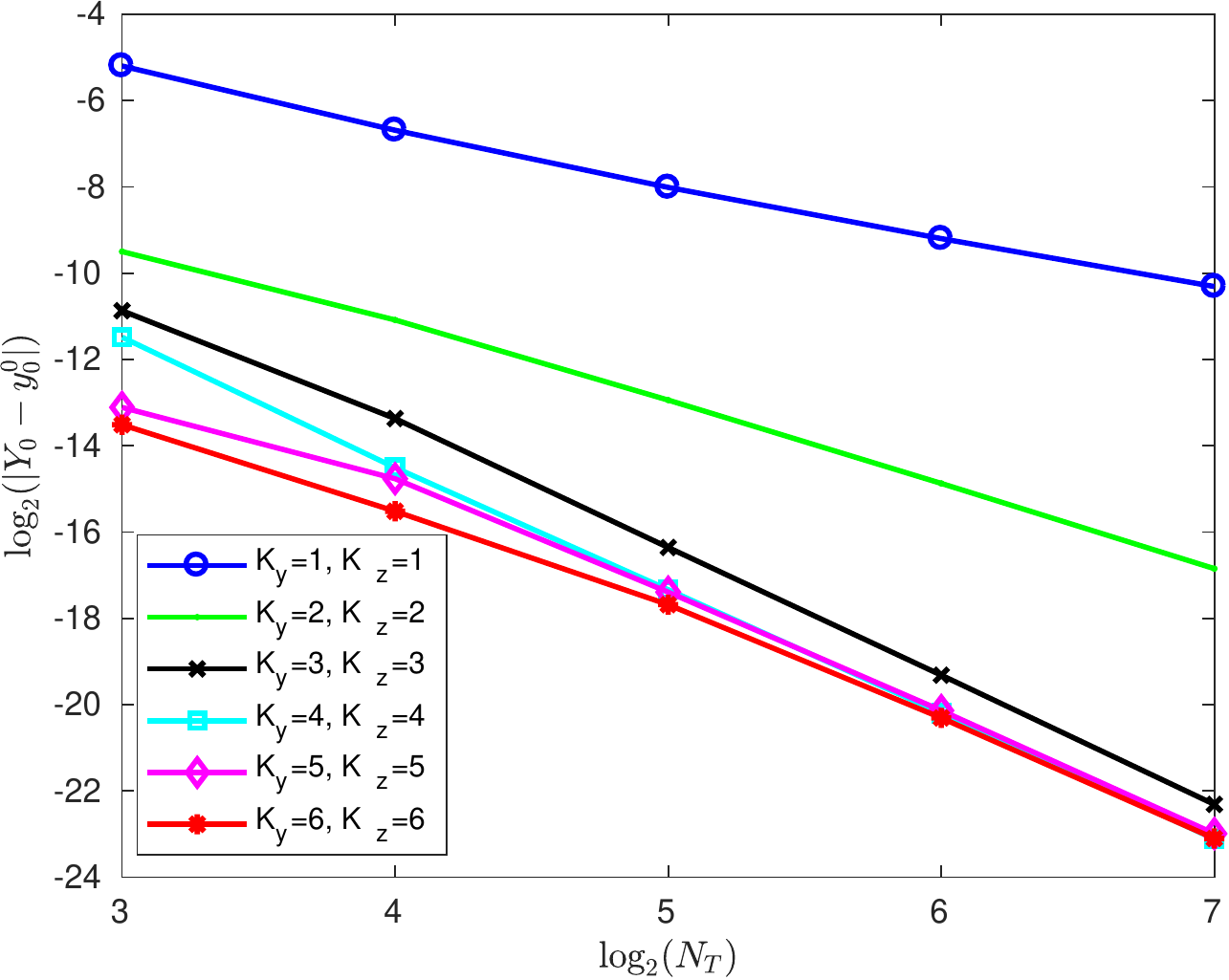}
 \subcaption{$Y$-component}
\end{subfigure}
 ~ 
 \begin{subfigure}[b]{0.47\textwidth}
 \includegraphics[width=\textwidth]{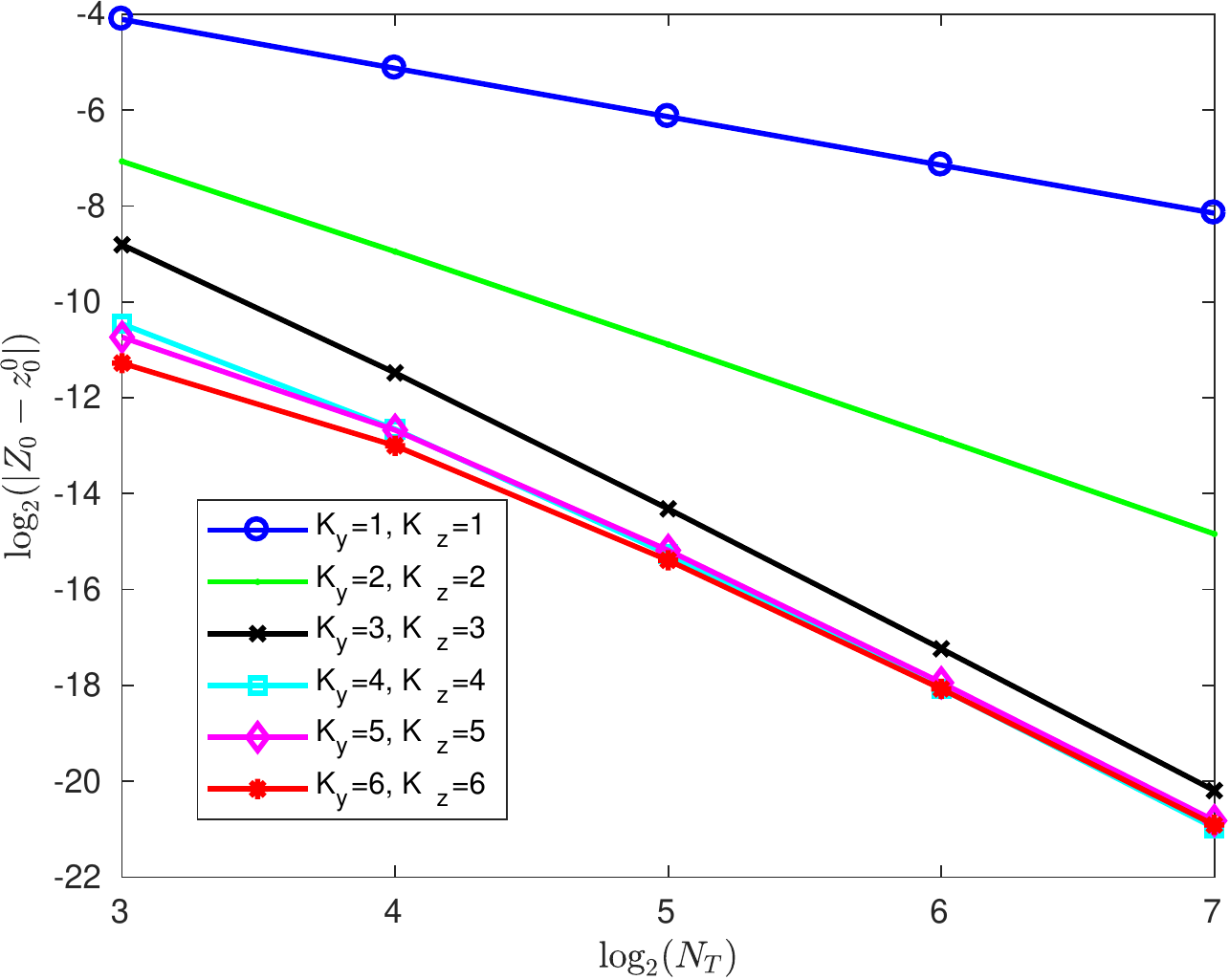}
 \subcaption{$Z$-component}
\end{subfigure}
 \caption{Plots of $\log_2\left(|Y_0-y_0^0|\right)$ and $\log_2\left(|Z_0-z_0^0|\right)$ with respect to $\log_2(N_T)$ for $K=1, \cdots 6$ for Example 2.}\label{fig:example_02}
\end{figure}

\begin{table}[h!]
\centering
\begin{tabular}{|c|c|c|c|c|c|c|}
\hline
\multicolumn{7}{|c|}{$|Z_0-z_0^0|$}  \\ \hline
 & $N_T=8$ & $N_T=16$ & $N_T=32$ & $N_T=64$ & $N_T=128$ & CR \\
\hline
$K_y=1, K_z=1, q=1$ & 5.80e-02 & 2.86e-02 & 1.42e-02 & 7.05e-03 & 3.52e-03 & 1.01 \\
\hline
$K_y=1, K_z=2, q=2$ & 9.45e-03 & 2.53e-03 & 6.54e-04 & 1.66e-04 & 4.20e-05 & 1.96 \\
\hline
$K_y=2, K_z=1, q=1$ & 5.99e-02 & 2.91e-02 & 1.43e-02 & 7.09e-03 & 3.53e-03 & 1.02 \\
\hline
$K_y=2, K_z=2, q=2$ & 7.45e-03 & 2.02e-03 & 5.28e-04 & 1.35e-04 & 3.41e-05 & 1.94 \\
\hline
$K_y=2, K_z=3, q=3$ & 2.25e-03 & 3.52e-04 & 4.91e-05 & 6.49e-06 & 8.35e-07 & 2.86 \\
\hline
$K_y=3, K_z=1, q=1$ & 5.99e-02 & 2.91e-02 & 1.43e-02 & 7.09e-03 & 3.53e-03 & 1.02 \\
\hline
$K_y=3, K_z=2, q=2$ & 7.46e-03 & 2.02e-03 & 5.28e-04 & 1.35e-04 & 3.41e-05 & 1.95 \\
\hline
$K_y=3, K_z=3, q=3$ & 2.23e-03 & 3.50e-04 & 4.90e-05 & 6.48e-06 & 8.34e-07 & 2.85 \\
\hline
$K_y=3, K_z=4, q=3$ & 6.84e-04 & 1.53e-04 & 2.53e-05 & 3.63e-06 & 4.86e-07 & 2.63 \\
\hline
$K_y=4, K_z=1, q=1$ & 5.99e-02 & 2.91e-02 & 1.43e-02 & 7.09e-03 & 3.53e-03 & 1.02 \\
\hline
$K_y=4, K_z=2, q=2$ & 7.44e-03 & 2.02e-03 & 5.28e-04 & 1.35e-04 & 3.41e-05 & 1.94 \\
\hline
$K_y=4, K_z=3, q=3$ & 2.26e-03 & 3.52e-04 & 4.91e-05 & 6.49e-06 & 8.35e-07 & 2.86 \\
\hline
$K_y=4, K_z=4, q=3$ & 7.10e-04 & 1.55e-04 & 2.54e-05 & 3.64e-06 & 4.86e-07 & 2.64 \\
\hline
$K_y=4, K_z=5, q=3$ & 5.94e-04 & 1.53e-04 & 2.69e-05 & 3.97e-06 & 5.40e-07 & 2.55 \\
\hline
$K_y=5, K_z=5, q=3$ & 5.86e-04 & 1.53e-04 & 2.69e-05 & 3.97e-06 & 5.40e-07 & 2.54 \\
\hline
$K_y=6, K_z=6, q=3$ & 4.03e-04 & 1.22e-04 & 2.33e-05 & 3.63e-06 & 5.08e-07 & 2.43 \\
\hline
\end{tabular}
\caption{Errors and convergence rates for Example 2, $T=1$}\label{tab:exp02z}
\end{table}
The given numerical results show that the proposed multi-step scheme works also well for a general nonlinear BSDE and is
a highly effective and accurate. Similar to Example 1, from Table \ref{tab:exp02y}, \ref{tab:exp02z} we can also observe that the results can be improved by increasing the number of steps.
And the estimated convergences rate are mostly consistent with the theoretical convergence order. Moreover, we observe that all estimated convergence rates are around $2.5$ for $K \geq 5.$
The reason for this is that the approximations (when $K \geq 5$) are too precise with $N_T=8.$
For this case we need to consider a greater value for $N_T$ in order to obtain an estimated rate close to $3.$ The plots of $\log_2\left(|Y_0-y_0^0|\right)$ and $\log_2\left(|Z_0-z_0^0|\right)$
with respect to $\log_2(N_T)$ are displayed in Figure \ref{fig:example_02}.

\paragraph{The Black-Scholes model}
In this example we compute the price of a European call option $V(t, S_t)$ by a BSDE where the underlying asset follows a geometric Brownian motion
\begin{equation}
 dS_t=\mu S_t\,dt + \sigma S_t dW_t.
\end{equation}
We assume that the asset pays dividends with the rate $d.$ The corresponding BSDE for the price of option can be derived
by setting up a self-financing portfolio $Y_t,$ which consists of $\pi_t$ assets and $Y_t-\pi_t$ bonds with risk-free return rate $r,$
which reads \cite{Karoui1997b}
\begin{equation}\label{eq:BSBSDE}
 \left\{
 \begin{array}{l}
 dS_t=\mu S_t\,dt + \sigma S_t \,dW_t,\\
-dY_t = \left(-r Y_t - \frac{\mu-r+d}{\sigma} Z_t\right)\,dt - Z_t\,dW_t,\\
\quad Y_T=\xi=\max(S_T-K, 0).
 \end{array}\right.
 \end{equation}
$Y_t$ is the option value $V(t, S_t),$ $Z_t$ corresponds to the hedging strategy, $Z_t=\sigma S_t \pi_t.$
We see that $S_t$ in \eqref{eq:BSBSDE} is a forward process, this type of BSDEs is called (uncoupled) forward backward stochastic
differential equation (FBSDE). The exact solution of \eqref{eq:BSBSDE} is given by the Black-Scholes model \cite{black1973}.
For $K=S=100, r=10\%, \mu=0.2, d=0, \sigma=0.25, T=0.1$ \footnote{We take the parameter values which are used in \cite{Ruijter2015} for comparison purpose.}, one obtains the exact solution $(Y_0,Z_0)=\left(3.65997,14.14823\right).$ 
In our experiment, for each time step we generate the grid point for $S$ by using 
the analytic solution of the geometric Brownian motion
\begin{equation}
 S_{i+1}=S_i\exp\left(\left(\mu-\frac{\sigma^2}{2}\right)h + \sigma \Delta x \right).
\end{equation}
Generally, one can use, e.g., the Euler or the Milstein method to simulate the forward process when there is no analytic solution available.

Note that the error analysis for the proposed methods relies on the smoothness assumptions of the initial data. However, in European option pricing,
the payoff function exhibits discontinuities at the strike price, this leads to a maximal error in the region of at-the-money.
For this problem, the smooth technqiue proposed by Kreiss et al. in \cite{Kreiss1970} has been widely used.
To further reduce the error caused by the missing smoothness we can e.g., start the multi-step algorithm without the (smoothed) initial data.
More precisely, we firstly smooth the initial data at $T.$ As mentioned before, for a $K$-step scheme we need to start with $K=1$ and choose a extremely small time step $\Delta t$ to compute
$(y_{\Gamma}^{N_T-j}, z_{\Gamma}^{N_T-j})$ for $j=1,\cdots,K-1$ using the smoothed initial data. Then, for computing $(y_{\Gamma}^{N_T-K}, z_{\Gamma}^{N_T-K})$ we use $y_{\Gamma}^{N_T-j}$ and $z_{\Gamma}^{N_T-j}$
only for $j=1, 2, \cdots,  K-1$ (without $j=0,$ namely without initial data), this computation is done by a $(K-1)$-step scheme. Finally, we can run the $K$-step scheme
to compute $(y_{\Gamma}^{N_T-K-1}, z_{\Gamma}^{N_T-K-1})$ based on $(y_{\Gamma}^{N_T-j}, z_{\Gamma}^{N_T-j}), j=1, 2, \cdots,  K,$ and so on backwards until the initial time.
%Of course, when the time step size $h$ already sufficiently small.
%Therefore, as a preparation, we firstly smooth the initial data using the method proposed by Kreiss et al. in \cite{Kreiss1970}, with a sufficiently small grid spacing for
%a high rate of convergence. 
We report our numerical results in Table \ref{tal:bsy} and \ref{tab:bsz}.
\begin{table}[!h]
\centering
\begin{tabular}{|c|c|c|c|c|c|c|}
\hline
\multicolumn{7}{|c|}{$|Y_0-y_0^0|$}  \\ \hline
 & $N=8$ & $N=16$ & $N=32$ & $N=64$ & $N=128$ & CR \\
 \hline
$K_y=1, K_z=1, q=1$ & 6.35e-04 & 2.88e-04 & 1.33e-04 & 6.78e-05 & 3.36e-05 & 1.06 \\
\hline
$K_y=1, K_z=2, q=2$ & 8.63e-06 & 1.02e-06 & 3.83e-07 & 1.22e-07 & 2.46e-08 & 2.00 \\
\hline
$K_y=2, K_z=1, q=1$ & 3.73e-04 & 1.70e-04 & 7.61e-05 & 3.92e-05 & 1.95e-05 & 1.06 \\
\hline
$K_y=2, K_z=2, q=2$ & 4.83e-06 & 1.31e-06 & 3.13e-07 & 4.85e-08 & 2.13e-08 & 2.04 \\
\hline
$K_y=2, K_z=3, q=3$ & 4.52e-09 & 3.83e-09 & 5.38e-10 & 7.70e-11 & 1.16e-11 & 2.29 \\
\hline
$K_y=3, K_z=1, q=1$ & 3.11e-04 & 1.60e-04 & 7.89e-05 & 4.22e-05 & 2.15e-05 & 0.96 \\
\hline
$K_y=3, K_z=2, q=2$ & 4.08e-06 & 8.78e-07 & 2.34e-07 & 8.79e-08 & 1.27e-08 & 2.00 \\
\hline
$K_y=3, K_z=3, q=3$ & 2.43e-08 & 3.37e-09 & 4.23e-10 & 8.75e-11 & 7.13e-12 & 2.87 \\
\hline
$K_y=3, K_z=4, q=3$ & 2.38e-08 & 3.33e-09 & 4.18e-10 & 8.69e-11 & 7.11e-12 & 2.87 \\
\hline
$K_y=4, K_z=1, q=1$ & 2.30e-04 & 1.25e-04 & 5.25e-05 & 2.70e-05 & 1.32e-05 & 1.05 \\
\hline
$K_y=4, K_z=2, q=2$ & 2.70e-06 & 6.17e-07 & 2.36e-07 & 5.80e-08 & 1.50e-08 & 1.84 \\
\hline
$K_y=4, K_z=3, q=3$ & 1.04e-08 & 1.25e-09 & 3.00e-10 & 4.85e-11 & 4.80e-12 & 2.69 \\
\hline
$K_y=4, K_z=4, q=3$ & 1.01e-08 & 1.22e-09 & 2.95e-10 & 4.79e-11 & 4.78e-12 & 2.68 \\
\hline
$K_y=4, K_z=5, q=3$ & 1.01e-08 & 1.19e-09 & 2.92e-10 & 4.76e-11 & 4.77e-12 & 2.67 \\
\hline
$K_y=5, K_z=5, q=3$ & 9.36e-09 & 1.68e-09 & 2.76e-10 & 2.97e-11 & 4.60e-12 & 2.78 \\
\hline
$K_y=6, K_z=6, q=3$ & 2.85e-08 & 1.38e-09 & 3.14e-10 & 3.13e-11 & 2.12e-12 & 3.29 \\
\hline
\end{tabular}
\caption{Errors and convergence rates for the Black-Scholes model}\label{tal:bsy}
\end{table}

\begin{figure}[htbp!]
 \centering
 \begin{subfigure}[b]{0.44\textwidth}
 \includegraphics[width=\textwidth]{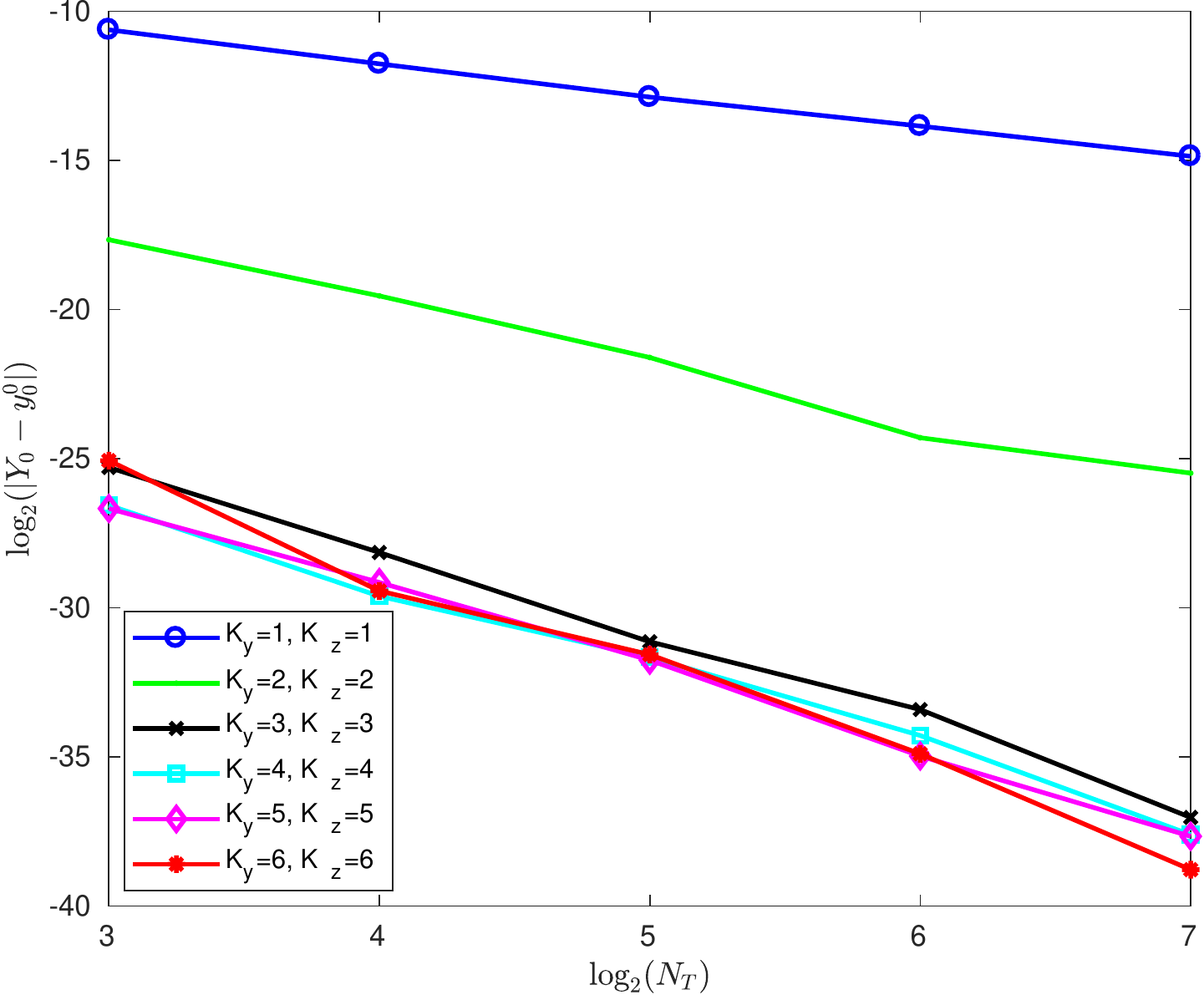}
 \subcaption{$Y$-component}
\end{subfigure}
 ~ 
 \begin{subfigure}[b]{0.47\textwidth}
 \includegraphics[width=\textwidth]{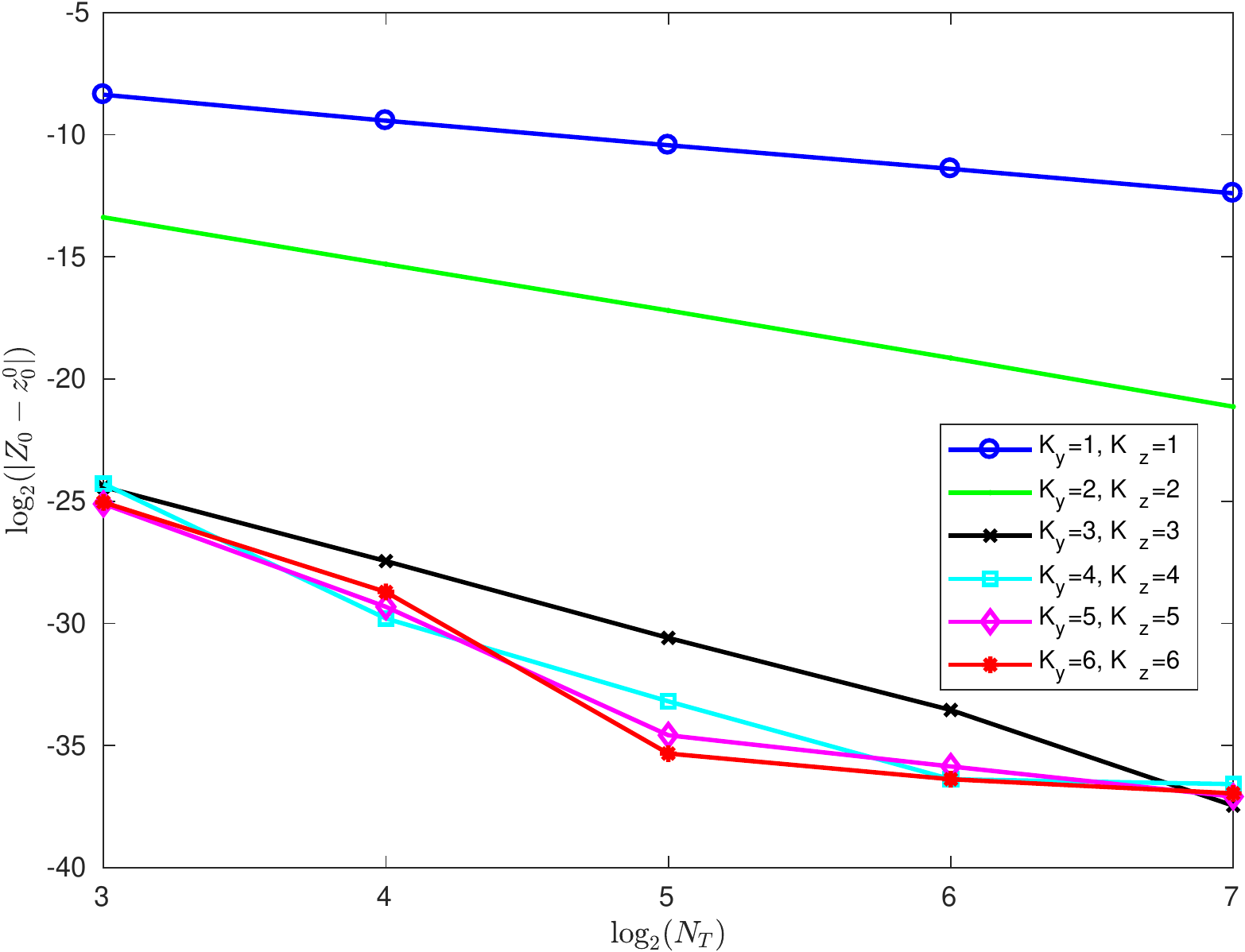}
 \subcaption{$Z$-component}
\end{subfigure}
 \caption{Plots of $\log_2\left(|Y_0-y_0^0|\right)$ and $\log_2\left(|Z_0-z_0^0|\right)$ with respect to $\log_2(N_T)$ for $K=1, \cdots 6$ for the example of the Black-Scholes model.}\label{fig:example_bs}
\end{figure}

\begin{table}[!h]
\centering
\begin{tabular}{|c|c|c|c|c|c|c|}
\hline
\multicolumn{7}{|c|}{$|Z_0-z_0^0|$}  \\ \hline
 & $N=8$ & $N=16$ & $N=32$ & $N=64$ & $N=128$ & CR \\
\hline
$K_y=1, K_z=1, q=1$ & 3.03e-03 & 1.45e-03 & 7.23e-04 & 3.70e-04 & 1.85e-04 & 1.00 \\
\hline
$K_y=1, K_z=2, q=2$ & 9.36e-05 & 2.46e-05 & 6.67e-06 & 1.73e-06 & 4.36e-07 & 1.93 \\
\hline
$K_y=2, K_z=1, q=1$ & 3.03e-03 & 1.46e-03 & 7.24e-04 & 3.71e-04 & 1.85e-04 & 1.00 \\
\hline
$K_y=2, K_z=2, q=2$ & 9.36e-05 & 2.48e-05 & 6.66e-06 & 1.73e-06 & 4.35e-07 & 1.93 \\
\hline
$K_y=2, K_z=3, q=3$ & 4.43e-08 & 5.05e-09 & 6.08e-10 & 7.92e-11 & 5.34e-12 & 3.20 \\
\hline
$K_y=3, K_z=1, q=1$ & 3.04e-03 & 1.46e-03 & 7.24e-04 & 3.71e-04 & 1.85e-04 & 1.00 \\
\hline
$K_y=3, K_z=2, q=2$ & 9.36e-05 & 2.48e-05 & 6.66e-06 & 1.73e-06 & 4.36e-07 & 1.93 \\
\hline
$K_y=3, K_z=3, q=3$ & 4.47e-08 & 5.45e-09 & 6.17e-10 & 7.98e-11 & 5.30e-12 & 3.22 \\
\hline
$K_y=3, K_z=4, q=3$ & 4.91e-08 & 9.42e-10 & 1.15e-10 & 1.08e-11 & 9.74e-12 & 3.10 \\
\hline
$K_y=4, K_z=1, q=1$ & 3.04e-03 & 1.46e-03 & 7.24e-04 & 3.71e-04 & 1.85e-04 & 1.00 \\
\hline
$K_y=4, K_z=2, q=2$ & 9.36e-05 & 2.48e-05 & 6.66e-06 & 1.73e-06 & 4.36e-07 & 1.93 \\
\hline
$K_y=4, K_z=3, q=3$ & 4.45e-08 & 5.42e-09 & 6.15e-10 & 7.93e-11 & 5.27e-12 & 3.22 \\
\hline
$K_y=4, K_z=4, q=3$ & 4.89e-08 & 1.07e-09 & 1.02e-10 & 1.12e-11 & 9.75e-12 & 3.12 \\
\hline
$K_y=4, K_z=5, q=3$ & 2.77e-08 & 1.09e-09 & 2.18e-11 & 1.65e-11 & 6.88e-12 & 3.05 \\
\hline
$K_y=5, K_z=5, q=3$ & 2.76e-08 & 1.49e-09 & 3.90e-11 & 1.61e-11 & 6.84e-12 & 3.05 \\
\hline
$K_y=6, K_z=6, q=3$ & 2.89e-08 & 2.27e-09 & 2.32e-11 & 1.12e-11 & 7.50e-12 & 3.15 \\
\hline
\end{tabular}
\caption{Errors and convergence rates for the Black-Scholes model}\label{tab:bsz}
\end{table}
From those tables, we clearly see that we have obtained surprisingly good accuracy.
The estimated convergence rates are again consistent with the theoretical convergence order.
Similar to the last two example, the approximation errors reduce mostly with the increasing
number of steps $K.$ We draw the plots of $\log_2\left(|Y_0-y_0^0|\right)$ and $\log_2\left(|Z_0-z_0^0|\right)$
with respect to $\log_2(N_T)$ in Figure \ref{fig:example_bs}.
\paragraph{Two-dimensional example}
For a two-dimensional example we consider the BSDE
\begin{equation*}
\left\{
\begin{array}{l}
-dY_t = \left(Y_t-\frac{Z^1_t}{2}-\frac{Z^2_t}{2}\right)\,dt - Z^1_t\,dW^1_t- Z^2_t\,dW^2_t,\\  
 \quad Y_T = \sin(W^1_T + W^2_T + T),
 \end{array}\right.
 \end{equation*}
with the analytic solution
 \begin{equation*}
 \left\{
 \begin{array}{l}
 Y_t = \sin(W^1_t + W^2_t + t),\\  
 Z_t = (\cos(W^1_t + W^2_t + t), \cos(W^1_t + W^2_t + t)),
 \end{array}\right.
 \end{equation*}
The exact solution of $(Y_0,Z^1_0,Z^2_0)$ is then $\left(0, 1, 1\right).$
The numerical approximations are reported in Table \ref{tab:twody} and \ref{tab:twodz}, which show that our multi-step scheme
is still quite highly accurate for solving a two-dimensional BSDE.
\begin{table}[!h]
\centering
\begin{tabular}{|c|c|c|c|c|c|c|}
\hline
\multicolumn{7}{|c|}{$|Y_0-y_0^0|$}  \\ \hline
  & $N=8$ & $N=16$ & $N=32$ & $N=64$ & $N=128$ & CR \\
\hline
$K_y=1, K_z=1, q=1$ & 1.32e-02 & 6.46e-03 & 3.18e-03 & 1.57e-03 & 7.81e-04 & 1.02 \\
\hline
$K_y=1, K_z=2, q=2$ & 4.72e-03 & 1.31e-03 & 3.45e-04 & 8.86e-05 & 2.24e-05 & 1.93 \\
\hline
$K_y=2, K_z=1, q=1$ & 1.22e-02 & 6.31e-03 & 3.17e-03 & 1.58e-03 & 7.88e-04 & 0.99 \\
\hline
$K_y=2, K_z=2, q=2$ & 1.83e-03 & 5.51e-04 & 1.48e-04 & 3.84e-05 & 9.82e-06 & 1.89 \\
\hline
$K_y=2, K_z=3, q=3$ & 3.97e-04 & 6.77e-05 & 9.74e-06 & 1.30e-06 & 1.65e-07 & 2.82 \\
\hline
$K_y=3, K_z=1, q=1$ & 8.59e-03 & 5.37e-03 & 2.94e-03 & 1.52e-03 & 7.76e-04 & 0.87 \\
\hline
$K_y=3, K_z=2, q=2$ & 1.48e-03 & 5.01e-04 & 1.42e-04 & 3.76e-05 & 9.69e-06 & 1.82 \\
\hline
$K_y=3, K_z=3, q=3$ & 3.94e-04 & 6.75e-05 & 9.72e-06 & 1.30e-06 & 1.64e-07 & 2.82 \\
\hline
$K_y=3, K_z=4, q=3$ & 1.88e-04 & 3.76e-05 & 5.68e-06 & 7.73e-07 & 9.78e-08 & 2.74 \\
\hline
$K_y=4, K_z=1, q=1$ & 5.44e-03 & 4.47e-03 & 2.70e-03 & 1.46e-03 & 7.61e-04 & 0.73 \\
\hline
$K_y=4, K_z=2, q=2$ & 1.14e-03 & 4.54e-04 & 1.36e-04 & 3.68e-05 & 9.60e-06 & 1.74 \\
\hline
$K_y=4, K_z=3, q=3$ & 2.91e-04 & 5.99e-05 & 9.21e-06 & 1.27e-06 & 1.63e-07 & 2.72 \\
\hline
$K_y=4, K_z=4, q=3$ & 1.90e-04 & 3.78e-05 & 5.69e-06 & 7.73e-07 & 9.77e-08 & 2.75 \\
\hline
$K_y=4, K_z=5, q=3$ & 1.42e-04 & 3.65e-05 & 5.99e-06 & 8.46e-07 & 1.09e-07 & 2.61 \\
\hline
$K_y=5, K_z=5, q=3$ & 1.39e-04 & 3.65e-05 & 5.99e-06 & 8.46e-07 & 1.09e-07 & 2.61 \\
\hline
$K_y=6, K_z=6, q=3$ & 8.12e-05 & 3.07e-05 & 5.49e-06 & 7.98e-07 & 1.05e-07 & 2.45 \\
\hline
\end{tabular}
\caption{Errors and convergence rates for the two-dimensional example}\label{tab:twody}
\end{table}
\begin{table}[!h]
\centering
\begin{tabular}{|c|c|c|c|c|c|c|}
\hline
\multicolumn{7}{|c|}{$\left(|Z^1_0-z_0^{0,1}|+|Z^2_0-z_0^{0,2}|\right)/2$}  \\ \hline
 & $N=8$ & $N=16$ & $N=32$ & $N=64$ & $N=128$ & CR \\
\hline
$K_y=1, K_z=1, q=1$ & 3.02e-02 & 4.77e-03 & 3.26e-03 & 1.87e-03 & 9.86e-04 & 1.12 \\
\hline
$K_y=1, K_z=2, q=2$ & 8.40e-03 & 2.31e-03 & 6.05e-04 & 1.54e-04 & 3.92e-05 & 1.94 \\
\hline
$K_y=2, K_z=1, q=1$ & 1.49e-02 & 3.92e-03 & 3.05e-03 & 1.82e-03 & 9.82e-04 & 0.90 \\
\hline
$K_y=2, K_z=2, q=2$ & 9.07e-03 & 2.51e-03 & 6.60e-04 & 1.69e-04 & 4.27e-05 & 1.94 \\
\hline
$K_y=2, K_z=3, q=3$ & 1.43e-03 & 2.08e-04 & 2.79e-05 & 3.59e-06 & 4.41e-07 & 2.92 \\
\hline
$K_y=3, K_z=1, q=1$ & 6.47e-03 & 2.99e-03 & 2.78e-03 & 1.75e-03 & 9.67e-04 & 0.63 \\
\hline
$K_y=3, K_z=2, q=2$ & 7.89e-03 & 2.37e-03 & 6.41e-04 & 1.66e-04 & 4.24e-05 & 1.89 \\
\hline
$K_y=3, K_z=3, q=3$ & 1.43e-03 & 2.08e-04 & 2.79e-05 & 3.59e-06 & 4.40e-07 & 2.92 \\
\hline
$K_y=3, K_z=4, q=3$ & 8.03e-04 & 1.23e-04 & 1.67e-05 & 2.16e-06 & 2.61e-07 & 2.90 \\
\hline
$K_y=4, K_z=1, q=1$ & 6.76e-03 & 2.13e-03 & 2.50e-03 & 1.68e-03 & 9.46e-04 & 0.60 \\
\hline
$K_y=4, K_z=2, q=2$ & 6.73e-03 & 2.21e-03 & 6.21e-04 & 1.64e-04 & 4.21e-05 & 1.84 \\
\hline
$K_y=4, K_z=3, q=3$ & 1.26e-03 & 1.98e-04 & 2.73e-05 & 3.55e-06 & 4.40e-07 & 2.88 \\
\hline
$K_y=4, K_z=4, q=3$ & 8.09e-04 & 1.23e-04 & 1.67e-05 & 2.16e-06 & 2.61e-07 & 2.90 \\
\hline
$K_y=4, K_z=5, q=3$ & 7.48e-04 & 1.30e-04 & 1.83e-05 & 2.41e-06 & 2.97e-07 & 2.84 \\
\hline
$K_y=5, K_z=5, q=3$ & 7.49e-04 & 1.30e-04 & 1.83e-05 & 2.41e-06 & 2.97e-07 & 2.84 \\
\hline
$K_y=6, K_z=6, q=3$ & 5.98e-04 & 1.18e-04 & 1.73e-05 & 2.31e-06 & 2.87e-07 & 2.77 \\
\hline
\end{tabular}
\caption{Errors and convergence rates for the two-dimensional example}\label{tab:twodz}
\end{table}

\begin{figure}[htbp!]
 \centering
 \begin{subfigure}[b]{0.47\textwidth}
 \includegraphics[width=\textwidth]{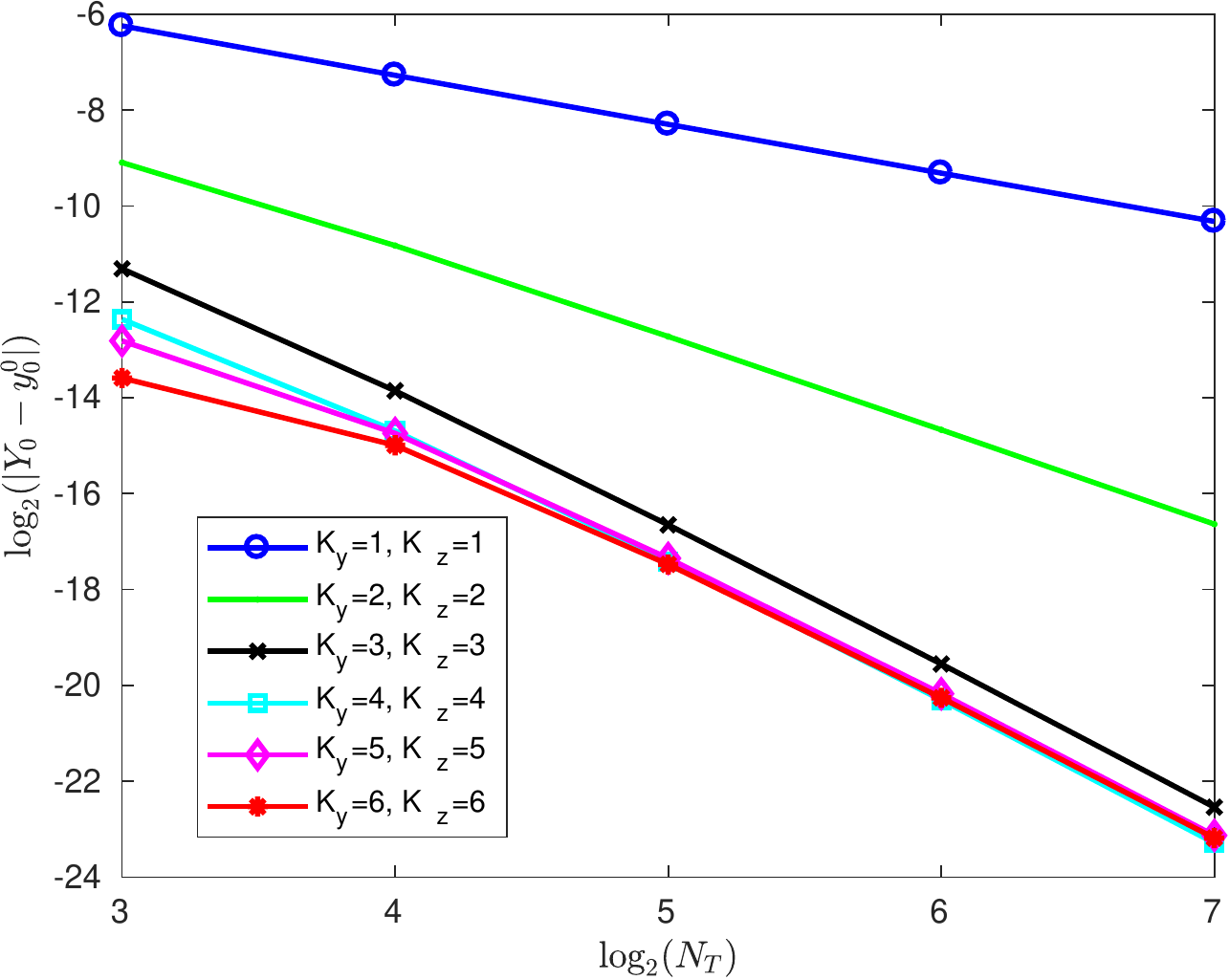}
 \subcaption{$Y$-component}
\end{subfigure}
 ~ 
 \begin{subfigure}[b]{0.47\textwidth}
 \includegraphics[width=\textwidth]{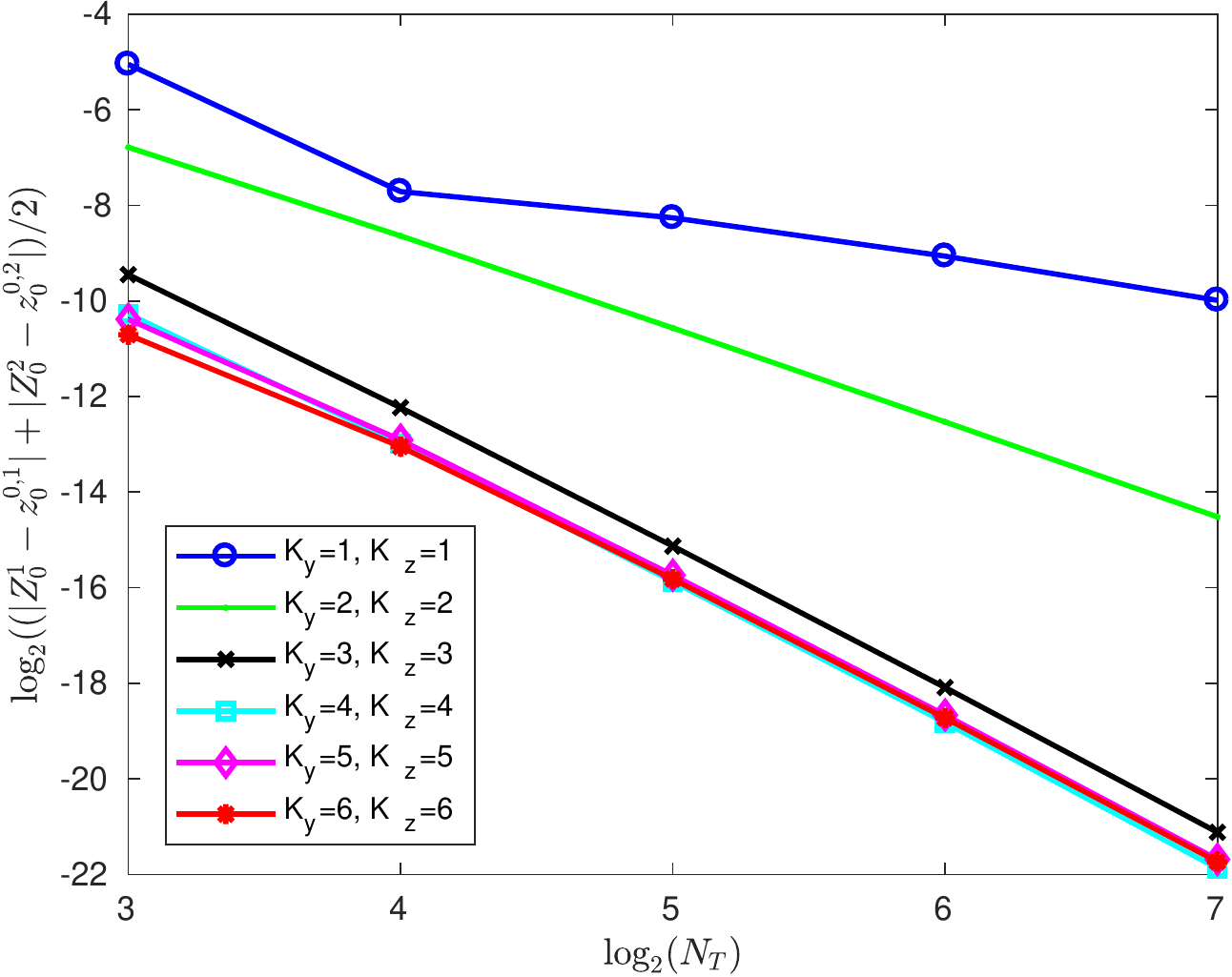}
 \subcaption{$Z$-component}
\end{subfigure}
 \caption{Plots of $\log_2\left(|Y_0-y_0^0|\right)$ and $\log_2\left((|Z^1_0-z_0^{0,1}|+|Z_0^2-z_0^{0,2}|)/2\right)$ with respect to $\log_2(N_T)$ for $K=1, \cdots 6$ for the two-dimensional example.}\label{fig:example_2d}
\end{figure}

As we have concluded for the one-dimensional examples above, in this two-dimensional example we see that
a smaller error value can be mostly achieved with a higher value of $K_y, K_z,$ i.e., more multi-steps.
The convergence rates are roughly consistent with the theoretical results in Theorem \ref{theorem:generalTheorem}.
The slight deviation comes from the quadratures and the two-dimensional interpolations.
The plots of $\log_2\left(|Y_0-y_0^0|\right)$ and $\log_2\left((|Z^1_0-z_0^{0,1}|+|Z_0^2-z_0^{0,2}|)/2\right)$
with respect to $\log_2(N_T)$ are given in Figure \ref{fig:example_2d}.
\section{Conclusion}
In this work, we adopt a multi-step scheme for solving BSDEs on time-space grids proposed in \cite{Zhao2010} by using the cubic spline interpolating
polynomials instead of the Lagrange interpolating polynomials in time. In \cite{Zhao2010} the number of multi-steps are limited, because the stability condition cannot be satisfied for a high number of time levels.
We find that our new proposed multi-step scheme allows for more multi-time-steps, which gives mostly a better approximation as our numerical results showed.
However, the convergence order of our scheme equals the one of scheme in \cite{Zhao2010}. The convergence order 
cannot be improved by using a higher value of $K.$ The reason for this is that a cubic spline is maximal fourth-order accurate.
Several numerical examples are provided to demonstrate the highly effectiveness and accuracy of our multi-step scheme for solving BSDEs.
In our proposed multi-step schemes, the computations among space grids at each time level are absolutly independent and should be thus parallelized.
Therefore, a GPU-based parallel computing is desirable for higher dimensional problems. This will be the task of future work.
\bibliographystyle{apalike}
\bibliography{mybibfile}

\begin{thebibliography}{}

\bibitem[Abramowitz and Stegun, 1972]{Abramowitz1972}
Abramowitz, M. and Stegun, I. (1972).
\newblock {\em Handbook of Mathematical Functions}.
\newblock Dover Publications.
\newblock Dover Books on Mathematics.

\bibitem[Bally, 1997]{Bally1997}
Bally, V. (1997).
\newblock Approximation scheme for solutions of bsde.
\newblock In Karoui, N.~E. and Mazliak, L., editors, {\em Backward stochastic
  differential equations}. Addison Wesley Longman, Harlow, UK.

\bibitem[Bender and Steiner, 2012]{Bender2012}
Bender, C. and Steiner, J. (2012).
\newblock Least-squares monte carlo for backward sdes.
\newblock {\em Numer. Methods Finance}, 12:257--289.

\bibitem[Bender and Zhang, 2008]{Bender2008}
Bender, C. and Zhang, J. (2008).
\newblock Time discretization and markovian iteration for coupled fbsdes.
\newblock {\em Ann. Appl. Probab.}, 18:143--177.

\bibitem[Black and Scholes, 1973]{black1973}
Black, F. and Scholes, M. (1973).
\newblock The pricing of options and corporate liabilities.
\newblock {\em J. Political Economy}, 81:637--654.

\bibitem[Bouchard and Touzi, 2004]{Bouchard2004}
Bouchard, B. and Touzi, N. (2004).
\newblock Discrete-time approximation and monte-carlo simulation of backward
  stochastic differential equations.
\newblock {\em Stoch. Proc. Appl.}, 111:175--206.

\bibitem[Crisan and Manolarakis, 2010]{Crisan2010}
Crisan, D. and Manolarakis, K. (2010).
\newblock Solving backward stochastic differential equations using the cubature
  method: Application to nonlinear pricing.
\newblock {\em SIAM J. FINAN. MATH}, 3(1):534--571.

\bibitem[Douglas et~al., 1996]{Douglas1996}
Douglas, J., Ma, J., and Protter, P. (1996).
\newblock Numerical methods for forward-backward stochastic differential
  equations.
\newblock {\em Ann. Appl. Probab.}, 6:940--968.

\bibitem[Gobet et~al., 2005]{Gobet2005}
Gobet, E., Lemor, J.~P., and Warin, X. (2005).
\newblock A regression-based monte carlo method to solve backward stochastic
  differential equations.
\newblock {\em Ann. Appl. Probab.}, 15:2172--2202.

\bibitem[Karoui et~al., 1997a]{Karoui1997a}
Karoui, N.~E., Kapoudjan, C., Pardoux, E., Peng, S., and Quenez, M.~C. (1997a).
\newblock Reflected solutions of backward stochastic differential equations and
  related obstacle problems for pdes.
\newblock {\em Ann. Probab.}, 25:702--737.

\bibitem[Karoui et~al., 1997b]{Karoui1997b}
Karoui, N.~E., Peng, S., and Quenez, M.~C. (1997b).
\newblock Backward stochastic differential equations in finance.
\newblock {\em Math. Finance}, 7(1):1--71.

\bibitem[Kreiss et~al., 1970]{Kreiss1970}
Kreiss, H.~O., Thom\'ee, V., and Widlund, O. (1970).
\newblock Smoothing of initial data and rates of convergence for parabolic
  difference equations.
\newblock {\em Commun. Pure Appl. Math}, 23(2):241--259.

\bibitem[Lemor et~al., 2006]{Lemor2006}
Lemor, J., Gobet, E., and Warin, X. (2006).
\newblock Rate of convergence of an empirical regression method for solving
  generalized backward stochastic differential equations.
\newblock {\em Bernoulli}, 12:889--916.

\bibitem[Lepeltier and Martin, 1997]{Lepeltier1997}
Lepeltier, J.~P. and Martin, J.~S. (1997).
\newblock Backward stochastic differential equations with continuous generator.
\newblock {\em Statist. Probab. Lett.}, 32(425--430).

\bibitem[Ma et~al., 2002]{Ma2002}
Ma, J., Protter, P., Mart\'in, J.~S., and Torres, S. (2002).
\newblock Numerical method for backward stochastic differential equations.
\newblock {\em Ann. Appl. Probab.}, 12:302--316.

\bibitem[Ma et~al., 1994]{Ma1994}
Ma, J., Protter, P., and Yong, J. (1994).
\newblock Solving forward-backward stochastic differential equations
  explicity-a four step scheme.
\newblock {\em Probab. Theory Related Fields}, 98(3):339--359.

\bibitem[Ma et~al., 2009]{Ma2009}
Ma, J., Shen, J., and Zhao, Y. (2009).
\newblock On numerical approximations of forward-backward stochastic
  differential equations.
\newblock {\em SIAM J. Numer. Anal.}, 46:2636--2661.

\bibitem[Ma and Zhang, 2005]{Ma2005}
Ma, J. and Zhang, J. (2005).
\newblock Representations and regularities for solutions to bsdes with
  reflections.
\newblock {\em Stoch. Proc. Appl.}, 115:539--569.

\bibitem[Milsetin and Tretyakov, 2006]{Milstein2006}
Milsetin, G.~N. and Tretyakov, M.~V. (2006).
\newblock Numerical algorithms for forward-backward stochastic differential
  equations.
\newblock {\em SIAM J. SCI. COMPUT.}, 28:561--582.

\bibitem[Pardoux and Peng, 1990]{Pardoux1990}
Pardoux, E. and Peng, S. (1990).
\newblock Adapted solution of a backward stochastic differential equations.
\newblock {\em System and Control Letters}, 14:55--61.

\bibitem[Pardoux and Peng, 1992]{Pardoux1992}
Pardoux, E. and Peng, S. (1992).
\newblock Backward stochastic differential equation and quasilinear parabolic
  partial differential equations.
\newblock {\em Lectures Notes in CSI.}, 176:200--217.

\bibitem[Peng, 1991]{Peng1991}
Peng, S. (1991).
\newblock Probabilistic interpretation for systems of quasilinear parabolic
  partial differential equations.
\newblock {\em Stochastics and Stochastic Reports}, 37(1--2):61--74.

\bibitem[Ruijter and Oosterlee, 2015]{Ruijter2015}
Ruijter, M.~J. and Oosterlee, C.~W. (2015).
\newblock A fourier cosine method for an efficient computation of solutions to
  bsdes.
\newblock {\em SIAM J. SCI. COMPUT.}, 37(2):A859--A889.

\bibitem[Teng, 2018]{Teng2018}
Teng, L. (2018).
\newblock A review of tree-based approaches to solve forward-backward
  stochastic differential equations.
\newblock {\em Preprint 18/03, University of Wuppertal}.

\bibitem[Zhang, 2004]{Zhang2004}
Zhang, J. (2004).
\newblock A numerical scheme for bsdes.
\newblock {\em Ann. Appl. Probab.}, 14:459--488.

\bibitem[Zhao et~al., 2006]{Zhao2006}
Zhao, W., Chen, L., and Peng, S. (2006).
\newblock A new kind of accurate numerical method for backward stochastic
  differential equations.
\newblock {\em SIAM J. SCI. COMPUT.}, 28(4):1563--1581.

\bibitem[Zhao et~al., 2009]{Zhao2009}
Zhao, W., Wang, J., and Peng, S. (2009).
\newblock Error estimates of the theta-scheme for backward stochastic
  differential equations.
\newblock {\em Discrete Contin. Dyn. Syst. Ser. B}, 12:905--924.

\bibitem[Zhao et~al., 2010]{Zhao2010}
Zhao, W., Zhang, G., and Ju, L. (2010).
\newblock A stable multistep scheme for solving backward stochastic
  differential equations.
\newblock {\em SIAM J. NUMER. ANAL.}, 48:1369--1394.

\end{thebibliography}
\end{document}